\documentclass[11pt]{article}

\usepackage[a4paper,margin=28mm]{geometry}
\usepackage[T1]{fontenc}
\usepackage{lmodern}
\usepackage{amsmath,amssymb,amsthm,mathtools,mathrsfs}
\usepackage{microtype}
\usepackage[colorlinks=true,linkcolor=blue,citecolor=blue,urlcolor=blue]{hyperref}
\hypersetup{
 pdftitle={Sharp Sobolev thresholds and rapid rotation for Navier--Stokes equations with fractional horizontal dissipation},
 pdfauthor={Zefu Feng, Yinghui Wang, Jiahong Wu},
 pdfsubject={Local well-posedness, conditional global existence, and Sobolev norm inflation},
 pdfkeywords={Navier--Stokes--Coriolis equations, horizontal dissipation, norm inflation, rapid rotation}
}

\numberwithin{equation}{section}
\newtheorem{theorem}{Theorem}[section]
\newtheorem{proposition}[theorem]{Proposition}
\newtheorem{lemma}[theorem]{Lemma}

\theoremstyle{remark}
\newtheorem{remark}[theorem]{Remark}
\allowdisplaybreaks[2]

\title{Sharp Sobolev thresholds and rapid rotation for Navier--Stokes
equations with fractional horizontal dissipation}
\author{Zefu Feng\thanks{School of Mathematical Sciences, Chongqing Normal
University, Chongqing 401331, China\protect\newline
(\texttt{zefufeng@mails.ccnu.edu.cn}).}
\and Yinghui Wang\thanks{MOE-LCSM, School of Mathematics and Statistics,
Hunan Normal University, Changsha 410081, China
(\texttt{yhwangmath@hunnu.edu.cn}).}
\and Jiahong Wu\thanks{Department of Mathematics, University of Notre Dame,
Notre Dame, IN 46556, USA (\texttt{jwu29@nd.edu}).}}
\date{}

\begin{document}
\maketitle

\begin{abstract}
We give a complete classification of the Sobolev spaces
$H^s(\mathbb R^3)$, $s>0$, in which the three-dimensional incompressible
Navier--Stokes equations with fractional horizontal dissipation of
order $2\alpha$, $0<\alpha\le1$, and a Coriolis force are well-posed.
The system has no vertical dissipation, and its scaling-critical index
is $s_c=5/2-2\alpha$.  For $0<\alpha<1$, we prove local Hadamard
well-posedness in $H^s$ for every $s\ge s_c$, including the critical
space itself; for $\alpha=1$, the admissible range is $s>1/2$.  For
all remaining positive Sobolev indices, we prove norm inflation for
classical solutions with smooth compactly supported initial data, at
each prescribed finite rotation rate.  The two results are
complementary and exhibit a dichotomy at the critical index:
$H^{s_c}$ is a well-posedness space when $0<\alpha<1$, whereas norm
inflation occurs in $H^{1/2}$ when $\alpha=1$.  Throughout the
well-posedness range, solutions are global for sufficiently small
initial data, with a smallness threshold independent of the rotation
rate, and for arbitrary fixed initial data at sufficiently large
rotation rates.  In the rapid-rotation limit, the difference from the
linear Stokes--Coriolis evolution converges to zero uniformly in time
and in the global dissipation norm; for spectrally localized data
separated from the vertical frequency axis, the rate is
$O(|\Omega|^{-\gamma})$ for every $0<\gamma<1/2$.  The well-posedness
theory is developed entirely in isotropic Sobolev spaces: for
$0<\alpha<1$, the energy estimate closes at the critical index without
loss of derivatives, by means of an isotropic Littlewood--Paley
decomposition, sharp trace inequalities, and frequency-localized
horizontal product estimates.  Below the critical index, norm
inflation follows from an Euler approximation.  At the endpoint
$(\alpha,s)=(1,1/2)$, the viscosity is invariant under the critical
rescaling and this perturbative argument is unavailable; we introduce
instead a construction based on a localized scalar shear and a
conormal analytic stability estimate that is uniform in the vertical
concentration parameters.
\end{abstract}



\section{Introduction}

We consider the Cauchy problem for the three-dimensional incompressible
Navier--Stokes equations with fractional horizontal dissipation and a
Coriolis force.  Horizontal viscosity with negligible vertical
viscosity arises in the study of rotating fluids; see
\cite{CheminDesjardinsGallagherGrenier2000,CheminDesjardinsGallagherGrenier2006}.
Since the dissipation acts only in the horizontal variables,
vertical transport has to be controlled without any vertical
smoothing, and the isotropic parabolic theory is not available.  This
paper determines exactly which Sobolev spaces $H^s(\mathbb R^3)$,
$s>0$, are well-posedness spaces for the resulting system, for the
whole range $0<\alpha\le1$ and at every finite rotation rate, and
shows how rapid rotation produces global solutions for arbitrarily
large initial data.

Let $x=(x_h,x_3)\in\mathbb R^2\times\mathbb R$, where $x_h=(x_1,x_2)$,
and set
\[
 H^s_\sigma(\mathbb R^3)
 =\{f\in H^s(\mathbb R^3;\mathbb R^3):\operatorname{div}f=0\},
 \qquad
 \widehat{\Lambda_h^{2\alpha}f}(\xi)
 =|\xi_h|^{2\alpha}\widehat f(\xi).
\]
The system reads
\begin{equation}\label{eq:NSC}
 \left\{
 \begin{aligned}
 &\partial_tu+\nu\Lambda_h^{2\alpha}u+\Omega e_3\times u
       +(u\cdot\nabla)u+\nabla p=0,
       &&(t,x)\in(0,\infty)\times\mathbb R^3,\\
 &\operatorname{div}u=0,\qquad u|_{t=0}=u_0,
 \end{aligned}
 \right.
\end{equation}
where $\nu>0$, $0<\alpha\le1$, $\Omega\in\mathbb R$, and $e_3=(0,0,1)$.
All data and solutions are real-valued.  The dissipative symbol
vanishes on the vertical frequency axis $\{\xi_h=0\}$.  Moreover,
$(e_3\times u,u)_{L^2}=0$, so the Coriolis force provides no direct
energy dissipation.  Its dispersive effect, combined with horizontal
heat damping away from that axis, is used to control the nonlinear
evolution for large $|\Omega|$.

When $\Omega=0$, the scaling
\[
 u(t,x)\longmapsto
 \lambda^{2\alpha-1}u(\lambda^{2\alpha}t,\lambda x)
\]
preserves the equation and the viscosity.  The corresponding
critical Sobolev index is
\[
                         s_c=\frac52-2\alpha.
\]
Our results give a complete answer in the scale of positive-order
Sobolev spaces.  The main contributions are the following.

\begin{itemize}
\item \emph{A sharp threshold that includes the critical space.}
For $0<\alpha<1$, we prove local Hadamard well-posedness in $H^s$ for
every $s\ge s_c$ and norm inflation for every $0<s<s_c$.  In
particular, the energy estimate closes in the isotropic critical
space $H^{s_c}$ without loss of derivatives, although the equation
has no vertical dissipation.  The earlier results for fractional
horizontal dissipation discussed below concern small data in $H^3$
or $H^2$.
\item \emph{A dichotomy at the critical index.}
For the horizontal Laplacian, $\alpha=1$, well-posedness holds for
$s>1/2$, while the critical space $H^{1/2}$ itself exhibits norm
inflation.  Thus the critical index is admissible precisely when
$0<\alpha<1$.  The endpoint $(\alpha,s)=(1,1/2)$ is not reached by
rescaling an Euler solution, because the viscosity is invariant under
the critical dilation; it requires a different mechanism, which we
construct from a localized scalar shear and a conormal analytic
stability estimate.
\item \emph{Global existence and stabilization by rotation.}
Throughout the well-posedness range, small data generate global
solutions with a smallness threshold independent of $\Omega$, and
every fixed datum, of arbitrary size, generates a global solution
once $|\Omega|$ is sufficiently large.  The nonlinear solution then
converges to the linear Stokes--Coriolis evolution uniformly on
$[0,\infty)$ and in the global dissipation norm, with an explicit
algebraic rate for spectrally localized data.
\end{itemize}

\subsection{Related work}

For the horizontal Laplacian, Chemin et~al.~\cite{CheminDesjardinsGallagherGrenier2000} established an anisotropic
Sobolev theory, and Iftimie \cite{Iftimie2002} subsequently extended
the uniqueness result.  Critical-space theories were developed by
Paicu \cite{Paicu2005} and Chemin and Zhang
\cite{CheminZhang2007}.  In these arguments, horizontal product
estimates and the identity
$\partial_3u_3=-\operatorname{div}_hu_h$ compensate for the lack of
vertical viscosity.  Paicu's critical existence theorem uses the
vertical Besov norm $\sum_j2^{j/2}\|\Delta_j^vf\|_2$, where
$\Delta_j^v$ denotes vertical frequency localization.  Burq et~al.~\cite{BurqLiuZhang2026} proved strong continuous dependence in
this critical Besov space.  The vertical dyadic $\ell^1$ summability
in this theory is stronger than the control supplied by the
isotropic $H^{1/2}$ norm; the endpoint norm inflation in
Theorem~\ref{thm:illposed-main} concerns the latter topology.

For fractional horizontal dissipation, Ji et~al.~\cite{JiLuoJiang2023} proved small-data $H^3$ stability for
$1/2<\alpha\le1$.  Bie et~al.~\cite{BieFangWangZhou2026} obtained small-data $H^3$ stability for
the directional operator $\Lambda_1^{2a}+\Lambda_2^{2b}$ throughout
$0<a,b\le1$, together with an $H^2$ result which, on the diagonal
$a=b=\alpha$, holds for $1/4<\alpha\le1$.  The equal-exponent
quadratic forms satisfy
\[
 c_\alpha\bigl(|\xi_1|^{2\alpha}+|\xi_2|^{2\alpha}\bigr)
 \le|\xi_h|^{2\alpha}
 \le C_\alpha\bigl(|\xi_1|^{2\alpha}+|\xi_2|^{2\alpha}\bigr).
\]
For the radial horizontal operator, Luo et~al.~\cite{LuoWuXu2026} proved small-data $H^3$ global stability for all
$0<\alpha\le1$ around a constant velocity, which can be removed by
a Galilean change of variables.  Theorem~\ref{thm:wellposed-main}
goes beyond these fixed integer orders: it gives local well-posedness
for arbitrary data and small-data global existence throughout the
sharp range \eqref{eq:wellposed-range}, which includes $H^2$ at the
borderline exponent $\alpha=1/4$ and the critical space $H^{s_c}$ for
every $0<\alpha<1$.

Time decay under classical horizontal dissipation was studied by Ji et~al.~\cite{JiWuYang2021}; Fujii
\cite{FujiiAnisotropic2023} analyzed the different asymptotic behavior
of the velocity components.  The fractional decay results in
\cite{BieFangWangZhou2026,LuoWuXu2026} require additional assumptions
on the low frequencies.  The limit in
\eqref{eq:linear-approximation} concerns rapid rotation with the
viscosity and initial datum fixed.  It gives approximation by the
linear Stokes--Coriolis evolution uniformly for all positive times
and in the global dissipation norm.

Rapid rotation provides another sufficient condition for global
existence.  In periodic domains, Babin et~al.~\cite{BabinMahalovNicolaenko1997} used averaging and resonance
analysis to obtain global regularity.  On $\mathbb R^3$, Chemin et~al.~\cite{CheminDesjardinsGallagherGrenier2000,CheminDesjardinsGallagherGrenier2002}
used dispersion to construct global solutions for a fixed large
initial datum and sufficiently fast rotation under classical
horizontal viscosity.  With zero two-dimensional background flow
and zero vertical viscosity, their linear approximation is stated
in anisotropic $H^{0,\sigma}$ spaces.
Theorem~\ref{thm:wellposed-main} treats the whole family of
fractional horizontal dissipations in isotropic Sobolev spaces, down
to the critical index when $0<\alpha<1$, and supplements the
qualitative limit \eqref{eq:linear-approximation} with the
quantitative estimate \eqref{eq:localized-stabilization-rate} for
spectrally localized data.  For rotating equations with full viscosity,
global well-posedness was studied in
\cite{HieberShibata2010,IwabuchiTakada2013,IwabuchiTakada2014}.
Related results include large-time behavior
\cite{EgashiraTakada2023}, large data in critical Fourier--Besov
spaces \cite{Fujii2024}, and isotropic fractional dissipation
\cite{AhnKimLee2022,SunDing2020}.

For the equations without rotation, Luo
\cite[Theorems~1.1 and~1.3]{Luo2026Sobolev} proved norm inflation for the
three-dimensional Euler equations when
$0<s<5/2$, and for the fully viscous Navier--Stokes equations when
$0<s<1/2$.  Both results use smooth, compactly supported,
divergence-free initial data and produce an arbitrarily large
$\dot H^s$ norm in an arbitrarily short time from arbitrarily small
$H^s$ data.  For the fully viscous equations, Luo
\cite[Theorem~1.1]{Luo2025Besov} also proved norm inflation in
supercritical Besov spaces $\dot B^s_{p,q}$ for
$s\ne0$, $1\le p,q\le\infty$, and $-3<s-3/p<-1$.
The critical space $H^{1/2}$ is outside this range.
With dissipation confined to the horizontal variables,
Theorem~\ref{thm:illposed-main} gives this form of norm inflation for
\eqref{eq:NSC} when $0<s<5/2-2\alpha$, for each prescribed finite
$\Omega$.  Its proof in this strict range uses
Luo's Euler norm-inflation theorem and treats horizontal dissipation
and rotation as perturbations after rescaling.  For $\alpha=1$,
our result also includes $s=1/2$, which lies outside the range of
Luo's Navier--Stokes theorem.  The viscosity remains fixed under
the critical rescaling, so the Euler approximation is no longer
available and the inflation has to be produced by the horizontally
viscous dynamics itself; the construction is described in
Section~\ref{subsec:outline}.  For $0<\alpha<1$, the critical index
$s=5/2-2\alpha$ belongs to the well-posedness range of
Theorem~\ref{thm:wellposed-main}.

\subsection{Main results}

Let $\mathbb P$ denote the Leray projector.  On divergence-free fields,
define the linear Stokes--Coriolis semigroup by
\begin{equation}\label{eq:S-definition}
 S_\Omega(t)
 =e^{-\nu t\Lambda_h^{2\alpha}}
       e^{-t\Omega\mathbb P(e_3\times\cdot)}.
\end{equation}
Theorem~\ref{thm:wellposed-main} gives local well-posedness,
small-data global existence, and global existence for a fixed datum
under sufficiently rapid rotation, together with the linear approximation.

\begin{theorem}[Local and conditional global well-posedness]
\label{thm:wellposed-main}
Let $\nu>0$ and suppose that
\begin{equation}\label{eq:wellposed-range}
 \begin{cases}
 0<\alpha<1,\quad s\ge \dfrac52-2\alpha,\\[2mm]
 \text{or}\\[-1mm]
 \alpha=1,\quad s>\dfrac12.
 \end{cases}
\end{equation}
Then the following assertions hold.

\begin{enumerate}
\renewcommand{\labelenumi}{\textup{(\roman{enumi})}}
\item
For every $\Omega\in\mathbb R$ and $u_0\in H^s_\sigma$, there is
$T>0$ and a unique distributional solution of \eqref{eq:NSC} in
the class
\begin{equation}\label{eq:solution-class}
 u\in C([0,T];H^s_\sigma),\qquad
 \Lambda_h^\alpha u\in L^2(0,T;H^s).
\end{equation}
For each reference datum, a common positive time can be chosen for
all sufficiently close $H^s_\sigma$ data.  The data-to-solution map
is continuous for the norm
\[
 \sup_{0\le t\le T}\|u(t)\|_{H^s}
       +\sqrt\nu\,\|\Lambda_h^\alpha u\|_{L^2(0,T;H^s)}.
\]
Uniqueness holds throughout every common interval in
\eqref{eq:solution-class}.

\item
There are constants $c_{s,\alpha},C_{s,\alpha}>0$, independent of
$\Omega$ and $\nu$, such that
$\|u_0\|_{H^s}\le c_{s,\alpha}\nu$ implies global existence and
\begin{equation*}
 \sup_{t\ge0}\|u(t)\|_{H^s}^2
 +\nu\int_0^\infty\|\Lambda_h^\alpha u(t)\|_{H^s}^2\,dt
 \le C_{s,\alpha}\|u_0\|_{H^s}^2.
\end{equation*}

\item
For each fixed $u_0\in H^s_\sigma$, there is
$\Omega_0=\Omega_0(u_0,\nu,\alpha,s)<\infty$ such that every
$|\Omega|\ge\Omega_0$ gives a unique global solution $u^\Omega$ with
\[
 u^\Omega\in C([0,\infty);H^s_\sigma),\qquad
 \sup_{t\ge0}\|u^\Omega(t)\|_{H^s}<\infty,\qquad
 \Lambda_h^\alpha u^\Omega\in L^2(0,\infty;H^s).
\]
Moreover,
\begin{equation}\label{eq:linear-approximation}
 \lim_{|\Omega|\to\infty}
 \left[
 \sup_{t\ge0}\|u^\Omega(t)-S_\Omega(t)u_0\|_{H^s}
 +\sqrt\nu\,\|\Lambda_h^\alpha
        (u^\Omega-S_\Omega(\cdot)u_0)\|_{L^2(0,\infty;H^s)}
 \right]=0.
\end{equation}

\item
If, in addition, there are $R>\kappa>0$ such that
\begin{equation*}
 \operatorname{supp}\widehat u_0
 \subset\{\xi:|\xi|\le R,\ |\xi_h|\ge\kappa\},
\end{equation*}
then for every $0<\gamma<1/2$ there are $C_\gamma>0$ and
$\Omega_\gamma\ge1$, depending on
$\gamma,s,\alpha,\nu,R,\kappa$ and $\|u_0\|_{H^s}$, such that
for $|\Omega|\ge\Omega_\gamma$,
\begin{equation}\label{eq:localized-stabilization-rate}
 \begin{split}
 &\sup_{t\ge0}\|u^\Omega(t)-S_\Omega(t)u_0\|_{H^s}\\
 &\qquad
 +\sqrt\nu\,\|\Lambda_h^\alpha
        (u^\Omega-S_\Omega(\cdot)u_0)\|_{L^2(0,\infty;H^s)}
 \le C_\gamma|\Omega|^{-\gamma}.
 \end{split}
\end{equation}
\end{enumerate}
\end{theorem}

The pressure gradient is recovered from the projected equation by
\[
 \nabla p=-(I-\mathbb P)
             \bigl(u\cdot\nabla u+\Omega e_3\times u\bigr).
\]
The global conclusions extend the local Hadamard solution; continuous
dependence holds on every finite interval of existence.  In
\eqref{eq:linear-approximation}, the linear and nonlinear solutions
have the same initial datum and rotation rate.

\begin{theorem}[Norm inflation in the complementary Sobolev range]
\label{thm:illposed-main}
Let $\nu>0$, $0<\alpha\le1$, and suppose that
\begin{equation*}
 0<s<\frac52-2\alpha,
 \qquad\text{or}\qquad
 (\alpha,s)=\left(1,\frac12\right).
\end{equation*}
For every fixed $\Omega\in\mathbb R$ and every
$\varepsilon,B,T_0>0$, there are
$u_0\in C_c^\infty(\mathbb R^3;\mathbb R^3)$ with
$\operatorname{div}u_0=0$ and $t_*\in(0,T_0)$ such that the
corresponding classical solution of \eqref{eq:NSC} exists through
$t_*$ and
\begin{equation}\label{eq:main-inflation}
             \|u_0\|_{H^s}<\varepsilon,
             \qquad
             \|u(t_*)\|_{\dot H^s}>B.
\end{equation}
Consequently, a solution map agreeing with classical solutions
cannot be continuous at zero from an $H^s_\sigma$ neighborhood
into $C([0,T];H^s_\sigma)$ on a common positive interval.
\end{theorem}

Theorems~\ref{thm:wellposed-main} and~\ref{thm:illposed-main} are
complementary: for each $0<\alpha\le1$, every $s>0$ is covered by
exactly one of them.  Hence the range \eqref{eq:wellposed-range} is
sharp in the scale of positive-order Sobolev spaces, at every fixed
finite rotation rate.

\begin{remark}[The integer Sobolev cases]
\label{rem:integer-cases}
Theorem~\ref{thm:wellposed-main} includes $H^3$ for every
$0<\alpha\le1$ and $H^2$ for $1/4\le\alpha\le1$.
For these integer orders, the equation also gives
$\partial_tu\in L^2_{\rm loc}(0,\infty;H^{s-1})$ for the global
solutions.  At $s=2$, Theorem~\ref{thm:illposed-main} proves norm
inflation for $0<\alpha<1/4$ at each fixed rotation rate.
Thus $\alpha=1/4$ is included in the sharp fixed-parameter
$H^2$ threshold.
\end{remark}

\begin{remark}[Rotation thresholds and critical topologies]
\label{rem:threshold-profile}
The threshold in Theorem~\ref{thm:wellposed-main}\textup{(iii)}
is selected after the datum and its spectral approximation have been
fixed, and can depend on the frequency distribution of the initial datum.
In Theorem~\ref{thm:illposed-main}, the finite rotation rate is fixed
first and the initial data producing norm inflation are then chosen.
At $\alpha=1$, well-posedness in a critical vertical Besov
space and conditional uniqueness in a critical energy class are
compatible with \eqref{eq:main-inflation}: smallness in the
isotropic $H^{1/2}$ norm does not provide the vertical dyadic
$\ell^1$ control used in that theory.
\end{remark}

\subsection{Outline of the proof}\label{subsec:outline}

The proof of Theorem~\ref{thm:wellposed-main} starts with nonlinear
estimates in an equivalent isotropic dyadic $H^s$ energy.  The
Coriolis term vanishes in each dyadic block.  For $0<\alpha<1$,
Sobolev trace inequalities and horizontal trilinear estimates give
the energy bound at $s=s_c$.  The case
$\alpha=1/2$ requires localized horizontal multiplier estimates.
In \eqref{eq:higher-sobolev-self}, the critical $H^{s_c}$ norm controls
the coefficient of the horizontal dissipation at every energy order.
In particular, the nonlinear energy estimate closes at $s=s_c$
without a loss of derivatives.
For $\alpha=1$, the vertical Sobolev embedding requires $s>1/2$.
Difference estimates one order below $H^s$, followed by a
frequency-tail estimate with slowly increasing weights, yield
uniqueness and continuous dependence in the solution class
\eqref{eq:solution-class}.

For global existence under rapid rotation, decompose a fixed datum
as $u_0=v_0+w_0$, where $w_0$ is small in $H^s$ and
$\operatorname{supp}\widehat v_0\subset
\{|\xi|\le R,\ |\xi_h|\ge\kappa\}$.
With $v=S_\Omega(t)v_0$ and $w=u-v$, the localized Coriolis
Strichartz estimates of Koh et~al.~\cite{KohLeeTakada2014a,KohLeeTakada2014b}, combined with horizontal
heat damping, give convergence to zero of $v$ in
$L^1(0,\infty;W^{k,\infty})$ and of $v\cdot\nabla v$ in
$L^1(0,\infty;H^s)$ as $|\Omega|\to\infty$, for a fixed integer
$k>s+3$.  At the energy level, the $w\cdot\nabla w$ term is absorbed by
the horizontal dissipation.  Estimating the forcing in $L^1$ in
time gives the global approximation
\eqref{eq:linear-approximation} and the rate
\eqref{eq:localized-stabilization-rate}.

For Theorem~\ref{thm:illposed-main} below $s_c$, the dilation
$u_N(t,x)=N^{3/2-s}V_N(N^{5/2-s}t,Nx)$ preserves the
$\dot H^s$ norm.  In the equation for $V_N$, the coefficients are
\[
 \nu_N=\nu N^{s+2\alpha-5/2},\qquad
 \Omega_N=\Omega N^{s-5/2}.
\]
Both tend to zero when $0<s<s_c$.  A classical Sobolev stability
estimate then transfers a fixed smooth Euler inflation solution to
\eqref{eq:NSC}.  At $(\alpha,s)=(1,1/2)$, one has $\nu_N=\nu$,
and this perturbative argument no longer applies.

At that endpoint, the inflation is produced within the
horizontally viscous dynamics.  The approximate solution is
constructed from a scalar horizontal parabolic equation, with a vertical mean
compensation that makes the velocity divergence-free and of finite
energy.  The vertical dyadic superposition yields an advected component with
small initial critical Sobolev norm and large amplitude near $x_3=0$.
The first parameter variation of the scalar advection--diffusion equation
satisfies \eqref{eq:endpoint-W-lower}, which yields the duality lower
bound \eqref{eq:endpoint-pairing-signal}.
A conormal analytic norm based on $x_3\partial_3$ gives
estimates for the Leray projector and Gaussian regularization at
a common analytic radius.  Horizontal localization makes the
vertical transport coefficient small, while horizontal dissipation
controls the loss of analytic radius in the error estimate.
Ordinary vertical derivatives satisfy a separate linear estimate,
which gives classical continuation for each fixed concentrating
profile.  The conormal analytic stability estimate in
Proposition~\ref{prop:endpoint-analytic-stability}, combined with
Lemma~\ref{lem:analytic-localization-approximation}, has constants uniform
in the vertical concentration parameters, so these parameters can be
chosen after the horizontal localization scales.  Duality against a uniformly bounded $\dot H^{-1/2}$
test function gives the lower bound for the Sobolev norm of the exact solution.

Section~\ref{sec:prelim} establishes the Sobolev energy and
difference estimates.  Section~\ref{sec:wellposedness} proves
Theorem~\ref{thm:wellposed-main}.  Sections
\ref{sec:subcritical-inflation} and \ref{sec:endpoint-inflation}
prove Theorem~\ref{thm:illposed-main} below $s_c$ and at
$(\alpha,s)=(1,1/2)$, respectively.

\section{Preliminaries and nonlinear estimates}\label{sec:prelim}

Let $\mathbb P$ denote the Leray projector,
\[
 \widehat{\mathbb Pf}(\xi)
  =\left(I-\frac{\xi\otimes\xi}{|\xi|^2}\right)\widehat f(\xi),
  \qquad \xi\ne0.
\]
Applying $\mathbb P$ to \eqref{eq:NSC} gives
\begin{equation}\label{eq:projected}
 \partial_tu+\nu\Lambda_h^{2\alpha}u
 +\Omega\mathbb P(e_3\times u)+\mathbb P(u\cdot\nabla u)=0.
\end{equation}
The scalar Fourier multipliers commute with $\mathbb P$.  For every
divergence-free $f$,
\begin{equation}\label{eq:coriolis-skew}
 \operatorname{Re}(\mathbb P(e_3\times f),f)_{L^2}=0.
\end{equation}
The identity remains valid after applying any scalar Fourier multiplier.

We use the unitary Fourier transform, with
$\langle\xi\rangle=(1+|\xi|^2)^{1/2}$, and write
\[
 \|f\|_{H^\sigma}=\|\langle\xi\rangle^\sigma\widehat f\|_2,
 \qquad \|f\|_{\dot H^\sigma}=\||\xi|^\sigma\widehat f\|_2.
\]
The same convention on $\mathbb R^2$ defines $H_h^\sigma$ and
$\dot H_h^\sigma$.  For a vector field, $f_h=(f_1,f_2)$;
$\mathcal F_hf$ denotes the Fourier transform in $x_h$ alone.
We set $D=-i\nabla$ and $\Lambda_h=(-\Delta_h)^{1/2}$;
$a(D)$ is the Fourier multiplier with symbol $a(\xi)$.
We write $L_3^qL_h^p$ for
$L^q(\mathbb R_{x_3};L^p(\mathbb R^2_{x_h}))$.
Horizontal homogeneous Sobolev norms are defined by the multiplier
$|\xi_h|^\gamma$; negative horizontal orders below are paired by
Fourier duality.

Fix a homogeneous isotropic Littlewood--Paley decomposition
$\{\Delta_j\}_{j\in\mathbb Z}$ with smooth real even radial cutoffs.
Write $f_j=\Delta_jf$ and $S_j=\sum_{k<j}\Delta_k$.
For vector fields, $f_{\ell,j}=\Delta_j(f_\ell)$ denotes the
$j$th block of the $\ell$th component.
Horizontal homogeneous blocks are denoted by $\Delta_\ell^h$.
The associated inhomogeneous isotropic blocks are
\begin{equation}\label{eq:dyadic-cutoffs}
 \Delta_{-1}^{\mathrm i}=\sum_{j<0}\Delta_j,
 \qquad \Delta_j^{\mathrm i}=\Delta_j\quad(j\ge0).
\end{equation}
Thus $S_k=\sum_{-1\le j<k}\Delta_j^{\mathrm i}$ for $k\ge0$;
the two decompositions use the same low-frequency cutoff.
A tilde denotes a fixed finite sum of neighboring blocks or sequence entries.
For $\sigma>0$ set
\begin{align*}
 \mathcal E_\sigma(f)^2
 &=\|f\|_2^2+\sum_{j\in\mathbb Z}2^{2j\sigma}\|f_j\|_2^2,
 &\mathcal D_{\sigma,\alpha}(f)
 &=\mathcal E_\sigma(\Lambda_h^\alpha f),
 \\
 \mathcal T_\sigma(u)
 &=\sum_{j\in\mathbb Z}2^{2j\sigma}
           (\Delta_j(u\cdot\nabla u),u_j)_{L^2}.
\end{align*}
These energies are equivalent to the inhomogeneous Sobolev norm and
its horizontal dissipation norm.  The $L^2$ nonlinear pairing vanishes.
Together with \eqref{eq:coriolis-skew}, this gives
\begin{equation*}
 \frac12\frac{d}{dt}\mathcal E_s(u)^2
       +\nu\mathcal D_{s,\alpha}(u)^2=-\mathcal T_s(u).
\end{equation*}
Generic constants may change from line to line.  They are independent
of $\Omega$ unless dependence is explicitly indicated.  All dyadic
calculations can first be performed for smooth fields with finite
frequency truncations and then passed to the limit.

For a divergence-free field whose Fourier support avoids the origin,
define the horizontal vector potential by
\begin{equation}\label{eq:vertical-potential}
 Q_\ell(f)=\partial_\ell\Delta^{-1}f_3
              -\partial_3\Delta^{-1}f_\ell,
 \qquad \ell=1,2,\qquad f_3=\operatorname{div}_hQ(f).
\end{equation}
The inverse Laplacian is used only after localization away from the origin.  The
fixed lowest frequencies are estimated separately in the proof of
Lemma~\ref{lem:diff-Hr}.

\subsection{A Sobolev trace inequality}

\begin{lemma}\label{lem:sobolev-trace}
For every $\sigma>1/2$,
\begin{equation*}
 \sup_{z\in\mathbb R}
       \|f(\cdot,z)\|_{H_h^{\sigma-1/2}}
             \le C_\sigma\|f\|_{H^\sigma(\mathbb R^3)}.
\end{equation*}
\end{lemma}
\begin{proof}
Fourier inversion in the vertical variable and Cauchy--Schwarz give
\begin{align*}
 |\mathcal F_hf(\xi_h,z)|^2
 &\le C\left(\int_{\mathbb R}
       (1+|\xi_h|^2+\xi_3^2)^{-\sigma}\,d\xi_3\right)
       \left(\int_{\mathbb R}
       (1+|\xi_h|^2+\xi_3^2)^\sigma
                   |\widehat f(\xi_h,\xi_3)|^2\,d\xi_3\right)\\
 &=C_\sigma(1+|\xi_h|^2)^{1/2-\sigma}
       \int_{\mathbb R}(1+|\xi_h|^2+\xi_3^2)^\sigma
                   |\widehat f(\xi_h,\xi_3)|^2\,d\xi_3.
\end{align*}
Multiply by $(1+|\xi_h|^2)^{\sigma-1/2}$ and integrate in $\xi_h$.
The vertical integral is finite because $\sigma>1/2$.
\end{proof}

\subsection{\texorpdfstring{The range $0<\alpha<1/2$}{The range 0 < alpha < 1/2}}

Assume $0<\alpha<1/2$. Define
\[
 X_\alpha=L_3^\infty L_h^{1/\alpha},\qquad
 Y_\alpha=L_3^2L_h^{2/(1-\alpha)},\qquad
 Z_\alpha=L_3^2L_h^{2/(1+\alpha)}.
\]
Mixed H\"older and the two-dimensional homogeneous Sobolev embedding give
\begin{equation}\label{eq:low-mixed}
 \|fg\|_{Z_\alpha}\leq\|f\|_{X_\alpha}\|g\|_{Y_\alpha},\qquad
 |(F,G)|\leq\|F\|_{Z_\alpha}\|G\|_{Y_\alpha},\qquad
 \|g\|_{Y_\alpha}\leq C_\alpha\|\Lambda_h^\alpha g\|_2.
\end{equation}

Since $s_c=5/2-2\alpha$,
$s_c-1=3/2-2\alpha>1/2$, Lemma~\ref{lem:sobolev-trace}
and the two-dimensional Sobolev embedding give
\begin{equation}\label{eq:low-trace}
 \|f\|_{X_\alpha}\le C_\alpha\|f\|_{H^{s_c-1}},
 \qquad
 \|\nabla f\|_{X_\alpha}\le C_\alpha\|f\|_{H^{s_c}}.
\end{equation}
The horizontal regularity supplied by the trace inequality is $1-2\alpha$.

\begin{lemma}[Mixed dyadic transport estimate]\label{lem:low-dyadic}
Let $0<\alpha<1/2$ and $\sigma>0$. For smooth divergence-free $u$,
\begin{equation}\label{eq:low-self}
 \sum_{j\in\mathbb Z}2^{2j\sigma}
 \left|\bigl(\Delta_j(u\cdot\nabla u),\Delta_ju\bigr)\right|
 \leq C_{\sigma,\alpha}\|\nabla u\|_{X_\alpha}
             \sum_j2^{2j\sigma}\|\Lambda_h^\alpha\Delta_ju\|_2^2.
\end{equation}
The principal transport term is canceled in each dyadic block
before the absolute value is taken. In particular, for $s\ge s_c$,
\begin{equation}\label{eq:low-self-Hs}
 |\mathcal T_s(u)|\leq C_{s,\alpha}\mathcal E_s(u)
                                     \mathcal D_{s,\alpha}(u)^2.
\end{equation}
\end{lemma}

\begin{proof}
Write $u_k=\Delta_ku$ and apply Bony's decomposition. For low--high
interactions, $|k-j|\leq C_0$. Separate the principal term
$S_{j-1}u\cdot\nabla u_j$, whose pairing with $u_j$ vanishes. The remaining
commutator satisfies
\[
 \|[\Delta_j,S_{k-1}u\cdot\nabla]u_k\|_{Z_\alpha}
 \leq C2^{-j}\|\nabla S_{k-1}u\|_{X_\alpha}
                      \|\nabla u_k\|_{Y_\alpha}
 \leq C\|\nabla u\|_{X_\alpha}\|u_k\|_{Y_\alpha}.
\]
This follows from the first moment of the convolution kernel. The
correction $S_{k-1}u-S_{j-1}u$ contains only finitely many comparable
blocks and obeys the same estimate. High--low terms are bounded by
$C\|\nabla u\|_{X_\alpha}\|u_k\|_{Y_\alpha}\|u_j\|_{Y_\alpha}$.

For high--high terms, $k\geq j-C_0$. Move the output cutoff to $u_j$ and
use the differentiated high factor in $X_\alpha$:
\[
 \left|\bigl(\Delta_j(u_k\cdot\nabla\widetilde u_k),u_j\bigr)\right|
 \leq C\|\nabla u\|_{X_\alpha}
                         \|u_k\|_{Y_\alpha}\|u_j\|_{Y_\alpha}.
\]
Cutoff kernels are bounded on $X_\alpha$, including its vertical
$L^\infty$ exponent. With $a_j=2^{j\sigma}\|u_j\|_{Y_\alpha}$, the
resulting sum is bounded by
\[
 C\|\nabla u\|_{X_\alpha}
   \sum_j a_j\sum_{k\geq j-C_0}2^{-\sigma(k-j)}a_k
 \leq C_\sigma\|\nabla u\|_{X_\alpha}\sum_j a_j^2.
\]
The geometric kernel is summable because $\sigma>0$. Finite-neighbor
sums are handled by Cauchy--Schwarz. Estimate \eqref{eq:low-mixed}
proves \eqref{eq:low-self}; \eqref{eq:low-trace} gives
\eqref{eq:low-self-Hs}. The argument uses only convolution-kernel bounds on the mixed spaces.
\end{proof}

\subsection{\texorpdfstring{The case $\alpha=1/2$}{The case alpha = 1/2}}

For the horizontal estimates, we distinguish whether the derivative
falls on the coefficient or on a high-frequency factor. For a smooth scalar function on $\mathbb R^3$, put
\[
 M_0(f)=\|\nabla_h f\|_{L_3^\infty L_h^2},\qquad
 M_B(f)=\sum_{\ell\in\mathbb Z}2^\ell
                \|\Delta_\ell^h f\|_{L_3^\infty L_h^2}.
\]
The latter dominates $\|f\|_{L_3^\infty\dot B^1_{2,1}(\mathbb R^2)}$.
This auxiliary coefficient norm will be applied to individual
isotropic blocks and differentiated low-frequency sums.

\begin{lemma}[Horizontal estimates and localized traces at $\alpha=1/2$]
\label{lem:critical-half-tools}
For smooth functions $F,G,H$ and $\ell\in\{1,2\}$,
\begin{align}
 \left|\int_{\mathbb R^3}(\partial_\ell F)GH\,dx\right|
 &\le C M_0(F)\|\Lambda_h^{1/2}G\|_2
                         \|\Lambda_h^{1/2}H\|_2,
 \label{eq:half-derivative-coefficient}\\
 \left|\int_{\mathbb R^3}F(\partial_\ell G)H\,dx\right|
 &\le C M_B(F)\|\Lambda_h^{1/2}G\|_2
                         \|\Lambda_h^{1/2}H\|_2.
 \label{eq:half-derivative-high}
\end{align}
The estimates are unchanged by independent translations of the factors.
For homogeneous isotropic blocks $f_q=\Delta_qf$, one has
\begin{align}
 M_0(f)&\le C\|f\|_{\dot H^{3/2}},
 \label{eq:half-sharp-trace}\\
 M_B(f_q)&\le C2^{3q/2}\|f_q\|_2,
 \label{eq:half-block-Besov}\\
 M_B(\nabla S_kf)&\le C2^k\|f\|_{\dot H^{3/2}}.
 \label{eq:half-gradient-low-Besov}
\end{align}
For the inhomogeneous decomposition, with $k\ge0$, the same estimates
hold with inhomogeneous norms and the usual low-block convention.
In addition,
\begin{align}
 M_B(S_kb)&\le C2^k\|b\|_{H^{1/2}},
 \label{eq:half-low-b-Besov}\\
 M_B(\partial_3S_kb)&\le C2^{2k}\|b\|_{H^{1/2}}.
 \label{eq:half-low-db-Besov}
\end{align}
\end{lemma}

\begin{proof}
For \eqref{eq:half-derivative-coefficient}, use horizontal H\"older
with exponents $(2,4,4)$, the embedding
$\dot H_h^{1/2}\hookrightarrow L_h^4$, and vertical exponents
$(\infty,2,2)$.
For \eqref{eq:half-derivative-high}, use the multiplier estimate
\[
 \|FH\|_{\dot H_h^{1/2}}
 \le C\|F\|_{\dot B^1_{2,1}(\mathbb R^2)}\|H\|_{\dot H_h^{1/2}}
\]
and duality against $\partial_\ell G\in\dot H_h^{-1/2}$.
This multiplier estimate follows from horizontal Bony decomposition:
the low--high paraproduct uses
$\dot B^1_{2,1}(\mathbb R^2)\hookrightarrow L_h^\infty$; for the high--low
paraproduct one uses
$\bigl\|\sum_{\ell<j}\Delta_\ell^hH\bigr\|_\infty\lesssim2^{j/2}\|H\|_{\dot H_h^{1/2}}$;
the high--high remainder is summed with the positive output exponent
$1/2$. Summation in $\ell^1$, as in the definition of $\dot B^1_{2,1}(\mathbb R^2)$, gives the multiplier bound.

The sharp trace \eqref{eq:half-sharp-trace} follows by Fourier
Cauchy--Schwarz in $\xi_3$ and the identity
\[
 \int_{\mathbb R}(|\xi_h|^2+\xi_3^2)^{-3/2}\,d\xi_3
 =2|\xi_h|^{-2}.
\]
After multiplication by $|\xi_h|^2$ and integration in $\xi_h$, this
gives the stated bound, uniformly in the vertical slice.

For \eqref{eq:half-block-Besov}, vertical Bernstein and the isotropic
frequency support give
\begin{align*}
 M_B(f_q)
 &\le C2^{q/2}\sum_{\ell\le q+C_0}
                        2^\ell\|\Delta_\ell^h f_q\|_2\\
 &\le C2^{q/2}
       \left(\sum_{\ell\le q+C_0}2^{2\ell}\right)^{1/2}
       \left(\sum_\ell\|\Delta_\ell^h f_q\|_2^2\right)^{1/2}
 \le C2^{3q/2}\|f_q\|_2.
\end{align*}
The sum over low horizontal frequencies converges geometrically,
without a logarithmic loss. Applying this estimate to $\nabla f_q$
and summing $q<k$ yields
\[
 M_B(\nabla S_k f)
 \le C\sum_{q<k}2^q\bigl(2^{3q/2}\|f_q\|_2\bigr)
 \le C2^k\|f\|_{\dot H^{3/2}}.
\]
For $b_q=\Delta_q^{\mathrm i}b$ and $k\ge0$, separate the low ball:
\begin{align*}
 M_B(S_kb)
 &\le C\|b_{-1}\|_2
       +C\sum_{0\le q<k}2^q\bigl(2^{q/2}\|b_q\|_2\bigr)
 \le C2^k\|b\|_{H^{1/2}},\\
 M_B(\partial_3S_kb)
 &\le C\|b_{-1}\|_2
       +C\sum_{0\le q<k}2^{2q}\bigl(2^{q/2}\|b_q\|_2\bigr)
 \le C2^{2k}\|b\|_{H^{1/2}}.
\end{align*}
These are \eqref{eq:half-low-b-Besov}--\eqref{eq:half-low-db-Besov}.
\end{proof}

\begin{lemma}[Weighted nonlinear energy estimate for $\alpha=1/2$]
\label{lem:critical-half-self}
Let $\sigma>0$, and let $\omega_j\ge1$ be nondecreasing, with
$\omega_{j+1}/\omega_j\le2^\epsilon$ for a fixed $\epsilon>0$. For smooth
divergence-free $u$,
\begin{equation*}
 \sum_{j\in\mathbb Z}\omega_j^22^{2\sigma j}
       \left|(\Delta_j(u\cdot\nabla u),u_j)_{L^2}\right|
 \le C_{\sigma,\omega}\|u\|_{\dot H^{3/2}}
       \sum_j\omega_j^22^{2\sigma j}
                         \|\Lambda_h^{1/2}u_j\|_2^2.
\end{equation*}
The principal transport term is canceled before estimating each block pairing.
The constant depends only on $\sigma$ and the fixed adjacent-ratio
bound, not on the values of the weights. The coefficient is the unweighted critical norm for every
$\sigma>0$ and for bounded or unbounded admissible weight sequences.
\end{lemma}

\begin{proof}
Set $d_j=\|\Lambda_h^{1/2}u_j\|_2$. Use the isotropic Bony decomposition.
For the potential in \eqref{eq:vertical-potential}, annular localization gives
\[
 \|\Lambda_h^{1/2}Q(f_k)\|_2
       \le C2^{-k}\|\Lambda_h^{1/2}f_k\|_2.
\]

For low--high transport, the principal term cancels. A horizontal
component of the commutator has a kernel moment $2^{-j}$, a coefficient
$\nabla S_{k-1}u_h$, and a horizontal derivative on a high factor.
Apply \eqref{eq:half-derivative-high} and
\eqref{eq:half-gradient-low-Besov}; since $|j-k|\le C_0$, its bound is
$C\|u\|_{\dot H^{3/2}} d_kd_j$.
For the vertical component, every component of $\nabla u_3$ is a
horizontal derivative of a velocity component, using
$\partial_3u_3=-\operatorname{div}_hu_h$.
Apply \eqref{eq:half-derivative-coefficient},
\eqref{eq:half-sharp-trace}, and the factor $2^k$ from
$\partial_3u_k$. The kernel moment cancels this factor, giving the
same bound. Comparable-frequency transport corrections use
\eqref{eq:half-block-Besov} and the potential $Q$ and also give
$C\|u\|_{\dot H^{3/2}} d_kd_j$.

For a horizontal high--low term, the derivative falls on the low
coefficient. Apply \eqref{eq:half-derivative-coefficient} directly,
using $M_0(S_{k-1}u)\le C\|u\|_{\dot H^{3/2}}$. For a vertical high--low term, write
$u_{3,k}=\operatorname{div}_hQ(u_k)$. The coefficient
$\partial_3S_{k-1}u$ has $M_B$ norm at most $C2^k\|u\|_{\dot H^{3/2}}$, while $Q(u_k)$
supplies $2^{-k}d_k$. Estimate \eqref{eq:half-derivative-high}
again gives $C\|u\|_{\dot H^{3/2}} d_kd_j$.

For high--high interactions with output $j\le k+C_0$, move the output
cutoff onto $u_j$. By \eqref{eq:half-block-Besov}, its $M_B$ norm is
at most $C2^{3j/2}\|u_j\|_2\le C\|u\|_{\dot H^{3/2}}$.
A horizontal derivative on a high factor is therefore handled by
\eqref{eq:half-derivative-high}. For a vertical derivative on a high
factor, first write the other high vertical component as
$\operatorname{div}_hQ(u_k)$. The factor $2^{-k}$ from $Q$ cancels the
factor $2^k$ from the vertical derivative. Both cases are bounded by
$C\|u\|_{\dot H^{3/2}} d_k\widetilde d_k$.

Consequently, each block pairing is bounded by
\[
 C\|u\|_{\dot H^{3/2}}\left(d_j\widetilde d_j+
                   \sum_{k\ge j-C_0}d_k\widetilde d_k\right).
\]
The finite-neighbor sum uses the bounded adjacent ratios of $\omega$.
For the high--high sum, monotonicity and $\sigma>0$ give
\[
 \sum_{j\le k+C_0}\omega_j^22^{2\sigma j}
 \le C_{\sigma,\omega}\omega_k^22^{2\sigma k}.
\]
Discrete Cauchy--Schwarz completes the proof.
\end{proof}

\subsection{\texorpdfstring{The range $1/2<\alpha<1$}{The range 1/2 < alpha < 1}}

\begin{lemma}[Horizontal trilinear estimates]
\label{lem:positive-high-horizontal}
Let $1/2<\alpha<1$.  For smooth scalar functions on
$\mathbb R^2$ and $\ell\in\{1,2\}$,
\begin{align}
 \left|\int fgh\,dx_h\right|
 &\leq C_\alpha\|f\|_{\dot H_h^{2-2\alpha}}
       \|g\|_{\dot H_h^{\alpha-1}}
       \|h\|_{\dot H_h^\alpha},
 \label{eq:positive-high-product}\\
 \left|\int (\partial_\ell f)gh\,dx_h\right|
 &\leq C_\alpha\|f\|_{\dot H_h^{2-2\alpha}}
       \|g\|_{\dot H_h^\alpha}
       \|h\|_{\dot H_h^\alpha}.
 \label{eq:positive-high-derivative}
\end{align}
In particular, the first estimate applies with $g=\partial_\ell G$ and
$\|g\|_{\dot H_h^{\alpha-1}}\leq\|\Lambda_h^\alpha G\|_2$.
These estimates remain valid after independent translations of the three
factors.  On $\mathbb R^3$, they can be integrated in $x_3$ with vertical
exponents $(\infty,2,2)$.
\end{lemma}

\begin{proof}
The two-dimensional homogeneous product law is
\[
 \|FG\|_{\dot H_h^{a+b-1}}
 \leq C_{a,b}\|F\|_{\dot H_h^a}\|G\|_{\dot H_h^b},
 \qquad a,b<1,\quad a+b>0.
\]
For the low--high paraproduct, Bernstein's inequality gives
\[
 \|S_kF\|_\infty\lesssim
 \sum_{\ell<k}2^{\ell(1-a)}2^{\ell a}\|\Delta_\ell^hF\|_2.
\]
The symmetric paraproduct is identical with $a,b$ exchanged.
For the high--high remainder, summation over output frequencies below
the input frequency uses $a+b>0$.  Cauchy--Schwarz in the dyadic indices
then gives the displayed product law.
Apply it first with $(a,b)=((2-2\alpha),\alpha-1)$.  Here $a,b<1$, $a+b=1-\alpha>0$, and $a+b-1=-\alpha$.  Duality proves
\eqref{eq:positive-high-product}.  Apply the same product law to $gh$ with
$a=b=\alpha$, obtaining $gh\in\dot H_h^{2\alpha-1}$, and pair against
$\partial_\ell f\in\dot H_h^{(2-2\alpha)-1}
 =\dot H_h^{1-2\alpha}$.  This proves
\eqref{eq:positive-high-derivative}.
\end{proof}

Set
\[
 M_\alpha(u):=\|u\|_{L_3^\infty\dot H_h^{2-2\alpha}}.
\]

Lemma~\ref{lem:sobolev-trace}, with
$\sigma=s_c=5/2-2\alpha>1/2$, gives
\begin{equation}\label{eq:positive-high-trace}
 M_\alpha(u)\le C_\alpha\|u\|_{H^{s_c}}.
\end{equation}

\begin{lemma}[Nonlinear energy estimate for $1/2<\alpha<1$]
\label{lem:positive-high-self}
Let $1/2<\alpha<1$ and $s>0$.  For every smooth divergence-free $u$,
\begin{equation}\label{eq:positive-high-self-M}
 \sum_{j\in\mathbb Z}2^{2js}
    \left|\bigl(\Delta_j(u\cdot\nabla u),u_j\bigr)_{L^2}\right|
 \leq C_{s,\alpha}M_\alpha(u)
       \sum_{j\in\mathbb Z}2^{2js}
       \|\Lambda_h^\alpha u_j\|_2^2.
\end{equation}
Consequently, if $s\ge5/2-2\alpha$,
\begin{equation}\label{eq:positive-high-self-Hs}
 \sum_j2^{2js}
    \left|\bigl(\Delta_j(u\cdot\nabla u),u_j\bigr)_{L^2}\right|
 \leq C_{s,\alpha}\mathcal E_s(u)
                         \mathcal D_{s,\alpha}(u)^2.
\end{equation}
\end{lemma}

\begin{proof}
Write $d_j=\|\Lambda_h^\alpha u_j\|_2$ and
$\widetilde d_j=\sum_{|k-j|\leq C_0}d_k$, where a fixed $C_0$ is enlarged
as necessary below.  The divergence-free condition gives, pointwise in
Fourier space,
\[
 |\widehat u_3(\xi)|\leq\frac{|\xi_h|}{|\xi|}|\widehat u(\xi)|.
\]
Thus the two block estimates needed for vertical interactions are
\begin{equation}\label{eq:positive-high-vertical-block}
 \|u_{3,j}\|_{L_3^2\dot H_h^{\alpha-1}}
       \leq C2^{-j}d_j,
 \qquad
 \|\partial_3u_j\|_{L_3^2\dot H_h^\alpha}
       \leq C2^jd_j.
\end{equation}
In the first estimate, incompressibility supplies the horizontal
factor needed for the negative homogeneous exponent.

Smooth isotropic cutoff kernels are bounded on the mixed spaces used
here.  In particular,
\begin{equation}\label{eq:positive-high-low-factor}
 M_\alpha(S_ju)+M_\alpha(u_j)\leq C M_\alpha(u),
 \qquad
 M_\alpha(\nabla S_ju)\leq C2^jM_\alpha(u).
\end{equation}

Apply the isotropic Bony decomposition to each component of
$u\cdot\nabla u$.  There are low--high terms with input index $k$ satisfying
$|k-j|\leq C_0$, high--low terms with the same index relation, and
high--high terms with $k\geq j-C_0$.

\smallskip
\noindent\emph{Low--high terms and cancellation.}
Commuting $\Delta_j$ through the low factor gives
$[\Delta_j,S_{k-1}u]\cdot\nabla u_k$ and a principal transport term.
The latter is written as
$S_{j-1}u\cdot\nabla u_j$ plus terms with all relevant isotropic frequencies
comparable to $2^j$.  The first term has zero $L^2$ pairing with $u_j$.
The comparable-frequency terms satisfy the high--high bounds
established in the last part of this proof.

For a horizontal component of the commutator, use its kernel formula
\[
 [\Delta_j,a]g(x)
 =\int K_j(y)\bigl(a(x-y)-a(x)\bigr)g(x-y)\,dy,
\]
and write the coefficient difference as
$-\int_0^1y\cdot\nabla a(x-\theta y)\,d\theta$.
The kernel moment is $\int|y||K_j(y)|\,dy\leq C2^{-j}$.
Use \eqref{eq:positive-high-product} with the low factor
$\nabla S_{k-1}u_h$, the horizontally differentiated high factor
$\nabla_hu_k$, and $u_j$.  By
\eqref{eq:positive-high-low-factor}, the resulting bound is
\[
 C2^{k-j}M_\alpha(u)d_kd_j
 \leq C M_\alpha(u)d_kd_j.
\]
For the vertical component of the commutator, every component of
$\nabla u_3$ is a horizontal derivative of a velocity component, since
$\partial_3u_3=-\operatorname{div}_hu_h$.
Use \eqref{eq:positive-high-derivative}, put
$\partial_3u_k$ in $L_3^2\dot H_h^\alpha$, and use
\eqref{eq:positive-high-vertical-block}.  The factors $2^{-j}$ and $2^k$
again cancel up to a fixed constant.  This gives the same bound.

\smallskip
\noindent\emph{High--low terms.}
For $u_{h,k}\cdot\nabla_hS_{k-1}u$, apply
\eqref{eq:positive-high-derivative} with the differentiated low factor
$S_{k-1}u$; both high factors have horizontal regularity $\alpha$.
The bound is $CM_\alpha(u)d_kd_j$.
For $u_{3,k}\partial_3S_{k-1}u$, apply
\eqref{eq:positive-high-product} with low factor
$\partial_3S_{k-1}u$.  Its mixed norm is at most $C2^kM_\alpha(u)$,
whereas \eqref{eq:positive-high-vertical-block} supplies $2^{-k}d_k$
from $u_{3,k}$.  This again yields $CM_\alpha(u)d_kd_j$.
Moving the output cutoff to $u_j$ changes none of these estimates.

\smallskip
\noindent\emph{High--high terms.}
Here $k\geq j-C_0$ and both inputs have frequency comparable to $2^k$.
Move the output cutoff to the output factor and place that factor in
$L_3^\infty\dot H_h^{2-2\alpha}$, with norm at most $CM_\alpha(u)$.
For a horizontal derivative on an input, use
\eqref{eq:positive-high-product} directly.  For a vertical derivative,
use \eqref{eq:positive-high-vertical-block}: the $2^{-k}$ from the vertical
velocity component cancels the $2^k$ from the vertical derivative of the
other input.  Each pairing is therefore bounded by
\[
 C M_\alpha(u)d_k\widetilde d_k.
\]
The coefficient norm is placed on the low output in this estimate.

Combining the three types of terms gives
\[
 |(\Delta_j(u\cdot\nabla u),u_j)|
 \leq C M_\alpha(u)
 \left(d_j\widetilde d_j+
       \sum_{k\geq j-C_0}d_k\widetilde d_k\right),
\]
with the comparable-frequency transport terms included by
enlarging $C_0$.  The first term is summed by Cauchy--Schwarz.  For the
second, $s>0$ gives
\[
 \sum_{j\leq k+C_0}2^{2js}\leq C_s2^{2ks}.
\]
A further finite-neighbor Cauchy--Schwarz estimate proves
\eqref{eq:positive-high-self-M}.  Estimate
\eqref{eq:positive-high-trace} proves
\eqref{eq:positive-high-self-Hs}.
\end{proof}

\subsection{The horizontal Laplacian}

For $\alpha=1$ we use the two-dimensional Gagliardo--Nirenberg
inequality.  Set
\[
 X(u)=\|u\|_{L_3^\infty L_h^2},\qquad
 Y(u)=\|\nabla_hu\|_{L_3^\infty L_h^2}.
\]
For the isotropic blocks $u_j$, write
$a_j=\|u_j\|_2$ and $d_j=\|\nabla_hu_j\|_2$.  Then
\begin{equation}\label{eq:positive-high-GN}
 \|u\|_{L_3^\infty L_h^4}
       \le C X(u)^{1/2}Y(u)^{1/2},
 \qquad
 \|u_j\|_{L_3^2L_h^4}\le C a_j^{1/2}d_j^{1/2}.
\end{equation}

\begin{lemma}\label{lem:laplacian-self}
For every $\sigma>0$ and every smooth divergence-free $u$,
\begin{align}
 &\sum_j2^{2j\sigma}
     |(\Delta_j(u\cdot\nabla u),u_j)_{L^2}|
 \notag\\
 &\quad\le C_\sigma\left[
 X(u)^{1/2}Y(u)^{1/2}
       \mathcal E_\sigma(u)^{1/2}\mathcal D_{\sigma,1}(u)^{3/2}
 +Y(u)\mathcal E_\sigma(u)\mathcal D_{\sigma,1}(u)\right].
 \label{eq:laplacian-self-mixed}
\end{align}
The principal low-frequency transport is canceled inside each block.
\end{lemma}
\begin{proof}
Apply the isotropic Bony decomposition, canceling
$(S_{j-1}u\cdot\nabla u_j,u_j)=0$ in each block.  Besides
\eqref{eq:positive-high-GN}, annular Bernstein and incompressibility give,
for every $k\in\mathbb Z$,
\[
 \|\partial_3u_k\|_{L_3^2L_h^4}
       \le C2^k(a_kd_k)^{1/2},\qquad
 \|u_{3,k}\|_2\le C2^{-k}d_k.
\]
For $|j-k|\le C_0$, the vertical high--low pairing satisfies
\[
 \left|(\Delta_j(u_{3,k}\partial_3S_{k-1}u),u_j)\right|
 \le C X(u)^{1/2}Y(u)^{1/2}d_k(a_jd_j)^{1/2}.
\]
Indeed, horizontal H\"older uses exponents $(2,4,4)$, the low factor
has $L_3^\infty L_h^4$ norm at most
$C2^kX(u)^{1/2}Y(u)^{1/2}$, and the two frequency factors cancel.
Horizontal commutators have the same bound: their kernel moment
$2^{-j}$ cancels the $2^k$ in
$\|\nabla S_{k-1}u_h\|_{L_3^\infty L_h^4}$.
For vertical commutators, use
$\|\nabla S_{k-1}u_3\|_{L_3^\infty L_h^2}\le CY(u)$ and the preceding
bound for $\partial_3u_k$.  Thus
\[
 \left|([\Delta_j,S_{k-1}u_3]\partial_3u_k,u_j)\right|
       \le CY(u)(a_kd_ka_jd_j)^{1/2}.
\]
Horizontal high--low terms obey this bound because
$\|\nabla_hS_{k-1}u\|_{L_3^\infty L_h^2}\le CY(u)$.

For high--high terms, $|k-k'|\le C_0$ and $j\le k+C_0$.
Move the output cutoff onto $u_j$, whose $L_3^\infty L_h^4$ norm
is bounded by $CX(u)^{1/2}Y(u)^{1/2}$.  A horizontal derivative
is estimated by $d_{k'}$; a vertical derivative contributes $2^{k'}$,
which is canceled by $2^{-k}d_k$ from the vertical velocity.  Each pairing is bounded by
\[
 C X(u)^{1/2}Y(u)^{1/2}
 \bigl[d_k(a_{k'} d_{k'})^{1/2}
             +d_{k'}(a_kd_k)^{1/2}\bigr].
\]
The comparable-frequency transport corrections satisfy the same estimate.
Since $\sigma>0$,
$\sum_{j\le k+C_0}2^{2\sigma j}\le C_\sigma2^{2\sigma k}$.
For these and the finite-neighbor terms, apply Cauchy--Schwarz
in the neighboring indices and then discrete H\"older:
\begin{align*}
 \sum_j2^{2\sigma j}a_j^{1/2}d_j^{3/2}
 &\le
 \Bigl(\sum_j2^{2\sigma j}a_j^2\Bigr)^{1/4}
 \Bigl(\sum_j2^{2\sigma j}d_j^2\Bigr)^{3/4},\\
 \sum_j2^{2\sigma j}a_jd_j
 &\le
 \Bigl(\sum_j2^{2\sigma j}a_j^2\Bigr)^{1/2}
 \Bigl(\sum_j2^{2\sigma j}d_j^2\Bigr)^{1/2}.
\end{align*}
Fixed neighboring shifts obey the same bounds.  The two right sides
are at most $\mathcal E_\sigma(u)^{1/2}
\mathcal D_{\sigma,1}(u)^{3/2}$ and
$\mathcal E_\sigma(u)\mathcal D_{\sigma,1}(u)$, respectively.
This proves \eqref{eq:laplacian-self-mixed}.
\end{proof}

\begin{proposition}[Nonlinear Sobolev energy estimate]\label{prop:sobolev-self}
Suppose that $0<\alpha<1$ and $s\ge5/2-2\alpha$, or that
$\alpha=1$ and $s>1/2$.  Then every smooth divergence-free $u$ satisfies
\begin{equation*}
 |\mathcal T_s(u)|
       \le C_{s,\alpha}\mathcal E_s(u)
                         \mathcal D_{s,\alpha}(u)^2.
\end{equation*}
For $0<\alpha<1$ and every energy order $\sigma>0$, one also has
\begin{equation}\label{eq:higher-sobolev-self}
 \sum_j2^{2j\sigma}
       |(\Delta_j(u\cdot\nabla u),u_j)_{L^2}|
       \le C_{\sigma,\alpha}\|u\|_{H^{5/2-2\alpha}}
                                      \mathcal D_{\sigma,\alpha}(u)^2.
\end{equation}
\end{proposition}
\begin{proof}
For $0<\alpha<1/2$, combine Lemma~\ref{lem:low-dyadic} with
\eqref{eq:low-trace}.  For $1/2<\alpha<1$, use
Lemma~\ref{lem:positive-high-self} and
\eqref{eq:positive-high-trace}.  At $\alpha=1/2$,
Lemma~\ref{lem:critical-half-self} with unit weights applies.
These arguments also prove \eqref{eq:higher-sobolev-self}.
For $\alpha=1$, Lemma~\ref{lem:sobolev-trace} gives, for $s>1/2$,
\[
 X(u)\le C_s\mathcal E_s(u),\qquad
 Y(u)\le C_s\mathcal D_{s,1}(u).
\]
Insert these bounds into \eqref{eq:laplacian-self-mixed} with
$\sigma=s$.
\end{proof}

\subsection{Weighted dyadic estimates}

Let $\omega_j\ge1$ be nondecreasing, equal to one for $j\le0$, and
satisfy $\omega_{j+1}\le2^\epsilon\omega_j$ for a fixed
$\epsilon>0$.  Define
\[
 \mathcal E_{\sigma,\omega}(u)^2
 =\|u\|_2^2+\sum_j\omega_j^22^{2j\sigma}\|u_j\|_2^2,
 \qquad
 \mathcal D_{\sigma,\alpha,\omega}(u)
 =\mathcal E_{\sigma,\omega}(\Lambda_h^\alpha u).
\]
The sequence may increase to infinity arbitrarily slowly.

\begin{lemma}\label{lem:weighted-self}
For $0<\alpha<1$ and $\sigma>0$,
\begin{equation}\label{eq:weighted-self-subunit}
 \sum_j\omega_j^22^{2j\sigma}
       |(\Delta_j(u\cdot\nabla u),u_j)_{L^2}|
 \le C_{\sigma,\alpha,\epsilon}
       \|u\|_{H^{5/2-2\alpha}}
                         \mathcal D_{\sigma,\alpha,\omega}(u)^2.
\end{equation}
For $\alpha=1$ and $\sigma>0$,
\begin{align}
 &\sum_j\omega_j^22^{2j\sigma}
       |(\Delta_j(u\cdot\nabla u),u_j)_{L^2}|
 \notag\\
 &\quad\le C_{\sigma,\epsilon}\left[
 X(u)^{1/2}Y(u)^{1/2}
      \mathcal E_{\sigma,\omega}(u)^{1/2}
                         \mathcal D_{\sigma,1,\omega}(u)^{3/2}
 +Y(u)\mathcal E_{\sigma,\omega}(u)
                         \mathcal D_{\sigma,1,\omega}(u)\right]
 \notag\\
 &\quad\le \frac{\nu}{2}\mathcal D_{\sigma,1,\omega}(u)^2
   +C_{\sigma,\epsilon}
       \left(\nu^{-3}X(u)^2Y(u)^2+\nu^{-1}Y(u)^2\right)
                              \mathcal E_{\sigma,\omega}(u)^2.
 \label{eq:weighted-self-laplacian}
\end{align}
The constants depend on the fixed neighboring-ratio bound and not
on the values of the weights.
\end{lemma}
\begin{proof}
For comparable frequencies, the weight ratios are uniformly bounded.
For high--high interactions,
\[
 \sum_{j\le k+C_0}\omega_j^22^{2j\sigma}
          \le C_{\sigma,\epsilon}\omega_k^22^{2k\sigma}.
\]
The proofs of Lemmas~\ref{lem:low-dyadic},
\ref{lem:positive-high-self}, and \ref{lem:laplacian-self}
therefore apply with these weighted sums.  For $\alpha=1/2$, the weighted assertion is already included
in Lemma~\ref{lem:critical-half-self}.
For the horizontal Laplacian, the same discrete H\"older estimates
give the first inequality in \eqref{eq:weighted-self-laplacian};
Young's inequality gives the second.
\end{proof}

\subsection{A lower-order difference estimate}

We estimate the difference of two solutions in the Sobolev space
one derivative below the solution space.
Throughout this subsection,
\[
0<\alpha<1,\qquad s\ge s_c=\frac52-2\alpha,
\qquad r=s-1.
\]
In particular, $s>1/2$ and $s+r=2s-1>0$.

Use the inhomogeneous blocks in \eqref{eq:dyadic-cutoffs} and put
$\lambda_j=2^{\max(j,0)}$.  For every real $q$, define
\[
 \mathfrak E_q(f)^2
 :=\sum_{j\geq-1}\lambda_j^{2q}
                      \|\Delta_j^{\mathrm i}f\|_2^2,
 \qquad
 \mathfrak D_{q,\beta}(f)
 :=\mathfrak E_q(\Lambda_h^\beta f).
\]
Thus $\mathfrak E_q\sim\|\cdot\|_{H^q}$ for every $q\in\mathbb R$.
Here $\sim$ denotes equivalence of norms.
In this subsection, $f_j=\Delta_j^{\mathrm i}f$,
$f_{\ell,j}=\Delta_j^{\mathrm i}(f_\ell)$, and
$S_k=\sum_{-1\le j<k}\Delta_j^{\mathrm i}$; in particular,
$S_k=0$ for $k\le-1$.

\begin{lemma}[Difference estimate below the solution regularity]
\label{lem:diff-Hr}
Let $0<\alpha<1$ and $s\ge5/2-2\alpha$, and set $r=s-1$.
There is $C=C(s,\alpha)$ such that, for smooth divergence-free vector
fields $a,b$,
\begin{align}
 &\sum_{j\geq-1}\lambda_j^{2r}
 \left|
    \bigl(\Delta_j^{\mathrm i}(a\cdot\nabla b),b_j\bigr)
 \right|
 +\sum_{j\geq-1}\lambda_j^{2r}
 \left|
    \bigl(\Delta_j^{\mathrm i}(b\cdot\nabla a),b_j\bigr)
 \right|\notag\\
 &\qquad\leq
 C\mathfrak E_s(a)
       \bigl(\mathfrak E_r(b)^2+\mathfrak D_{r,\alpha}(b)^2\bigr)
 +C\mathfrak E_r(b)\mathfrak D_{s,\alpha}(a)
                         \mathfrak D_{r,\alpha}(b).
 \label{eq:diff-Hr-main}
\end{align}
In the first sum, the vanishing principal transport pairing is removed
inside each block before taking its absolute value.
\end{lemma}

\begin{proof}[Proof for $1/2<\alpha<1$]
Use $M_\alpha$ from \eqref{eq:positive-high-trace} and write
$d_j(f)=\|\Lambda_h^\alpha f_j\|_2$.
The horizontal trilinear estimates
\eqref{eq:positive-high-product}--\eqref{eq:positive-high-derivative}
allow a horizontal derivative on any of the three factors.
The coefficient is estimated in
$\dot H_h^{2-2\alpha}$, and the other two factors are estimated in $\dot H_h^\alpha$.
The vertical exponents are $(\infty,2,2)$.
The trace estimate gives
\begin{equation}\label{eq:diff-M-trace}
 M_\alpha(a)\leq C\mathfrak E_s(a),\qquad
 M_\alpha\bigl((1-\Delta)^{-1/2}b\bigr)
       \leq C\mathfrak E_r(b).
\end{equation}
For a block with $j\geq0$, Bernstein's inequality also gives
\begin{equation}\label{eq:diff-output-Bernstein}
 M_\alpha(b_j)\leq C2^{js_c}\|b_j\|_2.
\end{equation}
Its low-block version has $2^j$ replaced by $1$.

For the potential $Q$ defined in \eqref{eq:vertical-potential},
annular localization and incompressibility give
\begin{equation}\label{eq:diff-Q-identity}
 f_{3,k}=\operatorname{div}_hQ(f_k),\qquad
 \|\Lambda_h^\alpha Q_\ell(f_k)\|_2
       \leq C2^{-k}d_k(f),\quad k\ge0,\ \ell=1,2.
\end{equation}
For an unlocalized field, apply $Q_\ell$ to $(I-S_0)b$.
The multiplier $Q_\ell(I-S_0)$ has order $-1$ and is smooth at the
origin.  Hence
\begin{align*}
 M_\alpha\bigl(Q_\ell((I-S_0)b)\bigr)
       &\leq C\mathfrak E_r(b),\\
 M_\alpha(S_kb)&\leq C2^k\mathfrak E_r(b),\qquad
 M_\alpha(\partial_3S_kb)\leq C2^{2k}\mathfrak E_r(b),\quad k\geq0.
\end{align*}
The cutoff multipliers $S_k(1-\Delta)^{1/2}$ and
$\partial_3S_k(1-\Delta)^{1/2}$ have $L^1$ kernel bounds
$C2^k$ and $C2^{2k}$.  Apply them to $(1-\Delta)^{-1/2}b$
and use \eqref{eq:diff-M-trace}.  All the
operators commute with $\Lambda_h^{2-2\alpha}$ and are bounded on the mixed
space defining $M_\alpha$.

\smallskip
\noindent\emph{The fixed lowest frequencies.}
Write $a^{\rm lo}=S_0a$ and $b^{\rm lo}=S_0b$.
These fields have Fourier support in a fixed ball, to which no
inverse derivative is applied.  For every fixed integer $m$,
\[
 \|b^{\rm lo}\|_{W^{m,\infty}}
       \leq C_m\mathfrak E_r(b).
\]
Consequently, multiplication by $b^{\rm lo}$ is bounded on $H^r$, and
\[
 \|b^{\rm lo}\cdot\nabla a\|_{H^r}
 \leq C\mathfrak E_r(b)\|a\|_{H^{r+1}}
 =C\mathfrak E_r(b)\|a\|_{H^s}.
\]
The term $a\cdot\nabla b^{\rm lo}$ has the same bound.
Since $a^{\rm lo}$ is divergence-free, the usual block
commutator calculation gives, for any real $r$,
\[
 \begin{aligned}
 &\left|\sum_{j\ge-1}\lambda_j^{2r}
   \bigl(\Delta_j^{\mathrm i}(a^{\rm lo}\cdot\nabla b),b_j\bigr)_{L^2}\right|\\
 &\quad+\left|\sum_{j\ge-1}\lambda_j^{2r}
   \bigl(\Delta_j^{\mathrm i}(b\cdot\nabla a^{\rm lo}),b_j\bigr)_{L^2}\right|
 \leq C\mathfrak E_s(a)\mathfrak E_r(b)^2.
 \end{aligned}
\]
Indeed, on high blocks the kernel moment $2^{-j}$ cancels the derivative
of $b_j$; the remaining bounded-frequency interactions are estimated
directly.  The same argument bounds the sum of the absolute block
pairings after principal transport cancellation.
Thus all terms requiring the inverse derivative at the origin contribute
only the first, nondissipative term on the right of
\eqref{eq:diff-Hr-main}.  The rest of the argument concerns annular
frequencies $2^k\gtrsim1$.

\smallskip
\noindent\emph{Low-frequency $a$ transporting high-frequency $b$.}
In the Bony term $S_{k-1}a\cdot\nabla b_k$, commute the output block
through the coefficient.  The principal term
$S_{j-1}a\cdot\nabla b_j$ has zero pairing with $b_j$.
For a horizontal component, the commutator kernel gives $2^{-j}$,
$M_\alpha(\nabla S_{k-1}a)\leq C2^kM_\alpha(a)$, and the remaining horizontal
derivative is handled by the horizontal trilinear estimates.  Since
$|j-k|\leq C_0$, the bound is
\[
 C M_\alpha(a)d_k(b)d_j(b).
\]
For the vertical component, each derivative in $\nabla a_3$ is a
horizontal derivative of a velocity component:
$\partial_\ell a_3$ already is one when $\ell=1,2$, while
$\partial_3a_3=-\operatorname{div}_ha_h$.
The horizontal trilinear estimates are therefore applied to that derivative;
$\partial_3b_k$ costs $2^k$, canceling the kernel moment $2^{-j}$.
The correction $S_{k-1}a-S_{j-1}a$ contains only comparable blocks.
Its horizontal part is estimated directly, and its vertical part uses
\eqref{eq:diff-Q-identity}, whose $2^{-k}$ cancels $\partial_3b_k$.
The resulting weighted sum is bounded by
$CM_\alpha(a)\mathfrak D_{r,\alpha}(b)^2$.

\smallskip
\noindent\emph{High-frequency $a$ and low-frequency derivatives of $b$.}
For $a_{h,k}\cdot\nabla_hS_{k-1}b$, use the horizontal trilinear estimates
with coefficient $S_{k-1}b$ and
$M_\alpha(S_{k-1}b)\leq C2^k\mathfrak E_r(b)$.
For $a_{3,k}\partial_3S_{k-1}b$, use
$a_{3,k}=\operatorname{div}_hQ(a_k)$ together with
\eqref{eq:diff-Q-identity} and
$M_\alpha(\partial_3S_{k-1}b)\leq C2^{2k}\mathfrak E_r(b)$.
Both estimates give
\[
 C2^k\mathfrak E_r(b)d_k(a)d_j(b),\qquad |k-j|\leq C_0.
\]
Their weighted sum is at most
$C\mathfrak E_r(b)\mathfrak D_{r+1,\alpha}(a)
\mathfrak D_{r,\alpha}(b)$, and $r+1=s$.

\smallskip
\noindent\emph{The paraproducts in $b\cdot\nabla a$.}
For low $b$ and high $a$, a horizontal derivative gives the preceding
$2^k\mathfrak E_r(b)$ coefficient.  For the vertical term, use
$S_{k-1}(I-S_0)b_3=\operatorname{div}_hS_{k-1}Q((I-S_0)b)$.
Its potential has $M_\alpha$ norm bounded by
$C\mathfrak E_r(b)$ and $\partial_3a_k$ costs $2^k$.
The fixed low-frequency part of $b_3$ was treated separately above.
This again gives
$C\mathfrak E_r(b)\mathfrak D_{s,\alpha}(a)
\mathfrak D_{r,\alpha}(b)$.
For high $b$ and low $a$, horizontal derivatives give
$CM_\alpha(a)d_k(b)d_j(b)$.  In the vertical term,
$b_{3,k}=\operatorname{div}_hQ(b_k)$ supplies $2^{-k}$, which cancels
$M_\alpha(\partial_3S_{k-1}a)\leq C2^kM_\alpha(a)$.
This contributes $CM_\alpha(a)\mathfrak D_{r,\alpha}(b)^2$.

\smallskip
\noindent\emph{High--high interactions and negative $r$.}
The two input indices are comparable to $k$ and
$j\leq k+C_0$.  Using divergence freeness, move the divergence to the
output factor.  For a horizontal component, estimate the low-frequency output in the
coefficient norm of the horizontal trilinear estimate.  This gives
\[
 C M_\alpha(b_j)d_k(a)\widetilde d_k(b).
\]
An output cutoff applied a second time only changes the constant.
For a vertical component, write the high vertical velocity as
$\operatorname{div}_hQ(a_k)$ or
$\operatorname{div}_hQ(b_k)$.  The factor $2^{-k}$ from $Q$ is combined
with $2^j$ from the vertical derivative of the output.  It yields the
additional bounded factor $2^{j-k}$.

For $j,k\geq0$, set
\[
 e_j=2^{jr}\|b_j\|_2,\qquad
 U_k=2^{ks}d_k(a),\qquad V_k=2^{kr}\widetilde d_k(b).
\]
By \eqref{eq:diff-output-Bernstein}, the weighted high--high sum is
bounded by a constant times
\begin{align*}
 \sum_{k\geq0}|U_kV_k|
 \sum_{0\leq j\leq k+C_0}
 |e_j|\,2^{j(r+s_c)-k(s+r)}
 &\leq C\sum_k|U_kV_k|
 \sum_{j\leq k+C_0}|e_j|
       2^{-(k-j)(2s-1)}2^{-j(s-s_c)}\\
 &\leq C\|e\|_{\ell^2}\|U\|_{\ell^2}\|V\|_{\ell^2}.
\end{align*}
Here $2s-1>0$ makes the geometric kernel square summable, and
$2^{-j(s-s_c)}\leq1$.  The output block $j=-1$ is bounded in the
same way, with coefficient $2^{-k(s+r)}$.
The high--high interaction estimate remains valid for $r<0$.

Collecting all terms and using \eqref{eq:diff-M-trace} proves
\eqref{eq:diff-Hr-main}.
\end{proof}

\begin{proof}[Proof for $0<\alpha<1/2$]
In this case $s\ge5/2-2\alpha>3/2$ and $r=s-1>1/2$.  Use the mixed spaces $X_\alpha,Y_\alpha,Z_\alpha$ in
\eqref{eq:low-mixed}.
Horizontal Sobolev embedding and the trace inequality give
\[
 \|f\|_{Y_\alpha}\leq C\|\Lambda_h^\alpha f\|_2,
 \qquad
 \|b\|_{X_\alpha}\leq C\|b\|_{H^r},
 \qquad
 \|\nabla a\|_{X_\alpha}\leq C\|a\|_{H^s}.
\]
Lemma~\ref{lem:sobolev-trace} applies since
$r\ge3/2-2\alpha>1/2$; the horizontal embedding uses
$r-1/2\ge1-2\alpha$.

For $a\cdot\nabla b$, low--high transport cancellation and its kernel
commutator give
\[
 C\|\nabla a\|_{X_\alpha}\mathfrak D_{r,\alpha}(b)^2.
\]
For high $a$ and low $b$, use
$\|\nabla S_kb\|_{X_\alpha}\leq C2^k\|b\|_{X_\alpha}$; this gives
\[
 C\|b\|_{X_\alpha}\mathfrak D_{r+1,\alpha}(a)
                         \mathfrak D_{r,\alpha}(b).
\]
For the high--high remainder, move the divergence to the output.
The resulting summation kernel is $2^{-(k-j)(r+1)}$, which is summable.
For $b\cdot\nabla a$, the low $b$--high $a$ term has the same bound;
the remaining terms are bounded by
$C\|\nabla a\|_{X_\alpha}\mathfrak D_{r,\alpha}(b)^2$.
In the latter high--high sum the kernel is $2^{-(k-j)r}$, which is summable since $r>0$.
Mixed H\"older uses $X_\alpha\cdot Y_\alpha\subset Z_\alpha$ and the dual pairing of $Z_\alpha$
with $Y_\alpha$.  The fixed lowest frequencies are handled by the smooth-multiplication
argument in the $1/2<\alpha<1$ part of the proof of
Lemma~\ref{lem:diff-Hr}.  These estimates prove \eqref{eq:diff-Hr-main}.
\end{proof}

\begin{proof}[Proof for $\alpha=1/2$]
Here $s\ge3/2$ and $r=s-1\ge1/2$.  Write $d_k(f)=\|\Lambda_h^{1/2}f_k\|_2$.
The fixed lowest frequencies are estimated by smooth multiplication,
contributing $C\mathfrak E_s(a)\mathfrak E_r(b)^2$.
For the remaining part, the order of $Q(I-S_0)$ and the sharp
trace give
\begin{equation}\label{eq:half-Qb-trace}
 M_0\bigl(Q((I-S_0)b)\bigr)
 \le C\|Q((I-S_0)b)\|_{H^{3/2}}
 \le C\mathfrak E_r(b).
\end{equation}
This estimate controls the low vertical velocity coefficient.

Low $a$ transporting high $b$ is estimated as in the proof of
Lemma~\ref{lem:critical-half-self}: horizontal commutators use
$M_B(\nabla S_ka)\le C2^k\mathfrak E_s(a)$, and vertical commutators use a horizontal
derivative of the coefficient and $M_0(a)\le C\mathfrak E_s(a)$.
The weighted sum is $C\mathfrak E_s(a)\mathfrak D_{r,1/2}(b)^2$.

For high $a$ and low derivatives of $b$, horizontal interactions use
$M_B(S_kb)\le C2^k\mathfrak E_r(b)$ (or the derivative-on-coefficient estimate).
Vertical interactions use $a_{3,k}=\operatorname{div}_hQ(a_k)$,
$M_B(\partial_3S_kb)\le C2^{2k}\mathfrak E_r(b)$, and the factor $2^{-k}$ from $Q$.
Thus both types are bounded blockwise by
\[
 C2^k\mathfrak E_r(b)\,d_k(a)d_j(b),\qquad |j-k|\le C_0.
\]
The same bound holds for low $b_h$ acting on high $a$.
For the part of low $b_3$ outside the fixed lowest ball, use
$(I-S_0)b_3=\operatorname{div}_hQ((I-S_0)b)$ and
\eqref{eq:half-derivative-coefficient} with coefficient
$Q((I-S_0)b)$.
By \eqref{eq:half-Qb-trace}, its coefficient is $C\mathfrak E_r(b)$, and
$\partial_3a_k$ supplies $2^kd_k(a)$. This gives the same displayed bound.
Since $r+1=s$, all these terms sum to
$C\mathfrak E_r(b)\mathfrak D_{s,1/2}(a)\mathfrak D_{r,1/2}(b)$.

For high $b$ acting on low $a$, a horizontal derivative falls on the
coefficient and is controlled by $M_0(a)$. In the vertical term,
$Q(b_k)$ supplies $2^{-k}$ and
$M_B(\partial_3S_ka)\le C2^k\mathfrak E_s(a)$ supplies the reciprocal factor.
These terms contribute $C\mathfrak E_s(a)\mathfrak D_{r,1/2}(b)^2$.

Finally, consider high--high interactions at index $k$ and output
$j\le k+C_0$. The coefficient norm of the low-frequency output is bounded by
\[
 M_B(b_j)\le C\lambda_j^{3/2}\|b_j\|_2.
\]
Horizontal derivatives are treated by
\eqref{eq:half-derivative-high}. In vertical terms the factor $2^{-k}$
from the high vertical velocity potential cancels the $2^k$ derivative
on the other high factor. Equivalently, moving the vertical divergence
to the output gives an additional factor $2^{j-k}$.
For $j,k\ge0$, set
\[
 e_j=2^{rj}\|b_j\|_2,\qquad
 U_k=2^{sk}d_k(a),\qquad
 V_k=2^{rk}\widetilde d_k(b).
\]
The resulting coefficient satisfies, for $j,k\ge0$,
\[
 2^{j(r+3/2)-k(s+r)}
 =2^{-(k-j)(2s-1)}2^{-j(s-3/2)}
 \le 2^{-(k-j)(2s-1)}.
\]
Therefore
\[
 \sum_k|U_kV_k|\sum_{j\le k+C_0}
        2^{-(k-j)(2s-1)}|e_j|
 \le C\|e\|_{\ell^2}\|U\|_{\ell^2}\|V\|_{\ell^2}.
\]
The lowest output block has coefficient $2^{-k(2s-1)}$ and is included in
the same estimate. Collecting all interactions proves
\eqref{eq:diff-Hr-main}.
\end{proof}

\section{Local well-posedness and conditional global existence}
\label{sec:wellposedness}

Throughout this section, $(\alpha,s)$ satisfies the hypothesis of
Theorem~\ref{thm:wellposed-main}, and $r=s-1$.
All smallness conditions below are imposed on a common constant
$c_{s,\alpha}>0$, which may be reduced finitely many times.  Its final
value is used throughout the proof and in
Theorem~\ref{thm:wellposed-main}\textup{(ii)}.
Let $S_\Omega(t)$ denote the semigroup
of the linear part of \eqref{eq:projected}.  It preserves bounded Fourier
support and is contractive on every $H^\sigma$ on divergence-free fields.

\begin{proposition}[Local solutions]\label{prop:local-theory}
For every $u_0\in H^s_\sigma(\mathbb R^3)$ and every fixed $\Omega$,
there is $T>0$ and a unique distributional solution of \eqref{eq:NSC}
satisfying
\[
 u\in C([0,T];H^s_\sigma),\qquad
 \Lambda_h^\alpha u\in L^2(0,T;H^s).
\]
For data in a sufficiently small $H^s$ neighborhood of a reference
datum, $T$ can be chosen in common, and the solution map is continuous
into both displayed spaces.  Uniqueness holds on every common interval
in this solution class.
\end{proposition}

We first establish the lower-order stability used in the proof.
For integer $M>5/2$, the classical local construction follows from
Friedrichs approximation and the commutator estimate
\cite{KatoPonce1988}
\[
 \frac{d}{dt}\|u\|_{H^M}^2
       +2\nu\|\Lambda_h^\alpha u\|_{H^M}^2
       \le C_M\|\nabla u\|_\infty\|u\|_{H^M}^2.
\]
The Coriolis term has zero pairing by \eqref{eq:coriolis-skew}.
The corresponding classical solution continues while its $H^M$
norm remains bounded.  The construction below gives a common
existence interval at the lower regularity $s$.

\subsection{Uniqueness and stability with a common smooth background}

Fix a smooth divergence-free field $v$, and let $u_i=v+w_i$,
$i=1,2$, be two solutions of \eqref{eq:projected} satisfying
\[
 w_i\in C([0,T];H^s),\qquad
 \Lambda_h^\alpha w_i\in L^2(0,T;H^s),
 \qquad i=1,2.
\]
The difference $z=w_1-w_2$ solves
\[
 \partial_tz+\nu\Lambda_h^{2\alpha}z
 +\Omega\mathbb P(e_3\times z)
 +\mathbb P\bigl(w_1\cdot\nabla z+z\cdot\nabla w_2
                  +v\cdot\nabla z+z\cdot\nabla v\bigr)=0.
\]
Fix an integer $m>\max(s+1,|s-1|+2)$.
With $r=s-1$, the smooth-coefficient block estimate bounds the
two terms containing $v$ by
$C\|v(t)\|_{W^{m,\infty}}\mathfrak E_r(z)^2$.

For $0<\alpha<1$, Lemma~\ref{lem:diff-Hr} and Young's inequality imply
that there is $\eta=\eta(\nu,s,\alpha)>0$ such that, if
\[
 \sup_{0\leq t\leq T}
 \bigl(\|w_1(t)\|_{H^s}+\|w_2(t)\|_{H^s}\bigr)\leq\eta,
\]
then
\begin{align}
 \frac{d}{dt}\mathfrak E_r(z)^2
 +c\nu\mathfrak D_{r,\alpha}(z)^2
 \leq\Bigl[ C\|v(t)\|_{W^{m,\infty}}
 +C\sum_{i=1}^2\bigl(
    \mathfrak E_s(w_i)
    +\nu^{-1}\mathfrak D_{s,\alpha}(w_i)^2\bigr)\Bigr]
                    \mathfrak E_r(z)^2.
 \label{eq:diff-common-background}
\end{align}
The coefficient belongs to $L^1(0,T)$.
One may take $\eta=c_{s,\alpha}\nu$, independently of $\Omega$.

\paragraph{The endpoint $\alpha=1$.}
For a fixed $s>1/2$, choose $\beta\in(1/2,1)$ sufficiently close to $1$
that $s>5/2-2\beta$.  Apply Lemma~\ref{lem:diff-Hr} with $\beta$ to the
nonlinear terms, while retaining the actual dissipation of order $1$.
For every real $q$,
\[
 \mathfrak D_{q,\beta}(f)
 \leq C\bigl(\mathfrak E_q(f)+\mathfrak D_{q,1}(f)\bigr).
\]
Smallness of $\|w_i\|_{H^s}$ absorbs the resulting multiples of
$\mathfrak D_{r,1}(z)^2$.  Young's inequality treats the other mixed
terms.  Thus
\begin{align}
 \frac{d}{dt}\mathfrak E_r(z)^2
 +c\nu\mathfrak D_{r,1}(z)^2
 \leq C_{\nu,s}\Bigl[1+C\|v(t)\|_{W^{m,\infty}}
 +\sum_{i=1}^2\bigl(\mathfrak E_s(w_i)^2
                    +\mathfrak D_{s,1}(w_i)^2\bigr)\Bigr]
                       \mathfrak E_r(z)^2.
 \label{eq:diff-alpha1}
\end{align}
This coefficient is integrable on every finite interval.  The auxiliary
$\beta$ is only used in a nonlinear estimate; the equation still has
horizontal Laplacian dissipation.

In all cases, Gronwall gives uniqueness and the lower-order stability
estimate
\begin{equation}\label{eq:diff-stability}
 \|w_1-w_2\|_{L^\infty(0,T;H^{s-1})}
 +\sqrt\nu\|\Lambda_h^\alpha(w_1-w_2)
                       \|_{L^2(0,T;H^{s-1})}
 \leq C_T\|w_1(0)-w_2(0)\|_{H^{s-1}}.
\end{equation}
The constant is uniform on a set with fixed bounds for the unweighted
solution and dissipation norms and for $\int_0^T\|v(t)\|_{W^{m,\infty}}\,dt$.
For nonsmooth solutions, the estimate follows by summing the fixed-block
identities.  The estimates in Lemma~\ref{lem:diff-Hr} make the canceled
block pairings absolutely summable in $L^1(0,T)$, so this passage does
not presuppose a classical derivative of an $H^{s-1}$ solution.

\subsection{Construction and strong continuous dependence}

We construct the solution by approximating its small remainder and
use weighted energies to obtain convergence in the stated Sobolev topology.

\begin{proof}[Proof of Proposition~\ref{prop:local-theory}]
We first choose a common existence interval for the smooth approximations.
Fix an integer $M>\max(s,5/2)$ and a small adjacent-ratio exponent
$\epsilon>0$ for the frequency weights used below.
Here $M$ is the auxiliary Sobolev order used for classical continuation;
$m$ denotes the order in the smooth-coefficient estimates above.
Choose the remainder threshold to satisfy the difference estimates,
the weighted estimates, and the estimate at this fixed order $M$.
Choose a fixed smooth spectral part $v_0$ of $u_0$ such that
$w_0=u_0-v_0$ is sufficiently small in $H^s$, and set
$v=S_\Omega(t)v_0$.  Using the cutoffs $S_N$ fixed in
Section~\ref{sec:prelim}, let $u^N$ be the classical solution with
initial value $v_0+S_Nw_0$.
Write $w^N=u^N-v$.
Proposition~\ref{prop:sobolev-self} and the estimates for terms
involving the smooth field $v$ give
\begin{align*}
 \frac12\frac{d}{dt}\mathcal E_s(w^N)^2
 +(\nu-C\mathcal E_s(w^N))\mathcal D_{s,\alpha}(w^N)^2
 \leq C\|v(t)\|_{W^{m,\infty}}\mathcal E_s(w^N)^2
       +C\|v\cdot\nabla v\|_{H^s}\mathcal E_s(w^N).
\end{align*}
Choose $T>0$ so that the bound obtained by Gronwall's inequality keeps
$\mathcal E_s(w^N)$ below a fixed strict fraction of the smallness
threshold.  This choice is independent of $N$.

To justify existence up to that common time, use the fixed order $M$.
For $\alpha<1$, \eqref{eq:higher-sobolev-self} at order $M$
has coefficient $\|w^N\|_{H^{5/2-2\alpha}}$, bounded by the
small $H^s$ norm.  The threshold chosen above absorbs this contribution.
The $H^M$ energy is then bounded on
$[0,T]$, with a bound depending on $N$ but finite for every $N$.
For $\alpha=1$, the higher-order Gagliardo--Nirenberg estimate instead
gives, after Young's inequality, a coefficient bounded by
\[
 C\bigl(\nu^{-3}\|w^N\|_{H^s}^2+\nu^{-1}\bigr)
                  \|\nabla_hw^N\|_{H^s}^2.
\]
Its time integral is uniformly bounded by the low-order estimate.
Gronwall again controls $H^M$ on the common interval.
Classical continuation therefore shows that every $u^N$ exists on
$[0,T]$.  No regularity assumption above $H^s$ is imposed on $u_0$.

\paragraph{Slow frequency weights.}
Choose a nondecreasing sequence $\omega_j\geq1$, equal to $1$ for
$j\leq0$, with $\omega_j\to\infty$ and
$\omega_{j+1}\leq2^\epsilon\omega_j$ for a fixed small
$\epsilon>0$, so slowly that
\[
 \mathcal E_{s,\omega}(w_0)<\infty.
\]
Such a choice follows by choosing successive frequency thresholds at which the
unweighted squared tail is, for example, less than $2^{-4n}$ and then
increasing the weight in sufficiently long steps.
The increase may be smoothed to satisfy the fixed neighboring-ratio
bound.  A single sequence works for every sequence of initial data
converging in $H^s$, because its squared high-frequency tails tend to
zero uniformly.

Apply Lemma~\ref{lem:weighted-self} to control the weighted energy.
For $0<\alpha<1$, use \eqref{eq:weighted-self-subunit}.
For $\alpha=1$, use \eqref{eq:weighted-self-laplacian}; its
coefficient is integrable by the unweighted dissipation bound.
The common background $v$ has fixed bounded Fourier support.
At frequencies above this support, multiplication by $v$ couples only
comparable blocks, so the neighboring-ratio bound preserves the
smooth-coefficient estimate.  On the bounded-frequency region the
weights are uniformly bounded.  The forcing $v\cdot\nabla v$ also has bounded
Fourier support.  Consequently,
\begin{equation}\label{eq:construction-weighted-bound}
 \sup_N\left[
 \sup_{0\leq t\leq T}\mathcal E_{s,\omega}(w^N(t))^2
 +\nu\int_0^T\mathcal D_{s,\alpha,\omega}(w^N(t))^2\,dt
 \right]<\infty.
\end{equation}
For $\alpha<1$, the smallness threshold can be chosen uniformly for
all weights with the fixed neighboring-ratio bound.  At $\alpha=1$, weighted propagation follows by Gronwall.

\paragraph{Convergence in the solution space.}
The approximations solve the same equation, so
\eqref{eq:diff-stability} applies directly to $w^N-w^L$; there are no
truncated nonlinear residuals.  Since $S_Nw_0$ is Cauchy in $H^{s-1}$,
\[
 w^N\quad\hbox{is Cauchy in}\quad
 C([0,T];H^{s-1})
 \cap\{f:\Lambda_h^\alpha f\in L^2(0,T;H^{s-1})\}.
\]
For a fixed frequency level $K$, low frequencies satisfy
\[
 \|S_K(w^N-w^L)\|_{H^s}
 \leq C2^K\|w^N-w^L\|_{H^{s-1}},
\]
and the same inequality holds after applying $\Lambda_h^\alpha$.
By \eqref{eq:construction-weighted-bound}, the high-frequency tails in
$C_tH^s$ and in $L_t^2H^s$ after applying $\Lambda_h^\alpha$ are bounded
by $C\omega_K^{-1}$, uniformly in $N,L$.
First choose $K$ large, then choose $N,L$ large.  This proves convergence
in
\begin{equation*}
 C([0,T];H^s)
 \cap\{f:\Lambda_h^\alpha f\in L^2(0,T;H^s)\}.
\end{equation*}
In particular the limit has the required strong initial trace.
Strong $C_tL^2$ convergence gives
$u^N\otimes u^N\to u\otimes u$ in $C_tL^1$, which suffices for passage
to the equation in distributions.  The linear terms pass by the same
convergence and the dissipation bound.  Estimate \eqref{eq:diff-stability} identifies the limit uniquely
among solutions whose remainders satisfy the common smallness bound.

\paragraph{Continuous dependence.}
For $u_{0,n}\to u_0$ in $H^s$, use the same spectral part $v_0$ chosen
for $u_0$ and write $w_{0,n}=u_{0,n}-v_0$.
The strict smallness of the remainders gives a common interval and uniform
unweighted estimates for all sufficiently large $n$.
A common slowly increasing weight is available for $w_{0,n}$ and $w_0$.
Apply \eqref{eq:diff-stability} first and then the preceding
low-frequency/high-frequency decomposition.  One obtains
\[
 \|u_n-u\|_{C([0,T];H^s)}
 +\sqrt\nu\|\Lambda_h^\alpha(u_n-u)\|_{L^2(0,T;H^s)}
 \longrightarrow0.
\]
This proves strong continuous dependence in the solution class.

\paragraph{Uniqueness in the full solution class.}
Let two solutions in \eqref{eq:solution-class} have the same initial value.
Choose the same smooth spectral part of that value.
Strong $H^s$ continuity makes both remainders smaller than the
threshold in \eqref{eq:diff-common-background}, or in
\eqref{eq:diff-alpha1}, on a common short interval.
Estimate \eqref{eq:diff-stability} gives equality there.
At any time where they agree, repeat this argument with a smooth
spectral approximation of their common value.  The set of agreement
times is closed by continuity and extends to the right at each of its
points before the endpoint.  This proves uniqueness on the entire
common interval and completes the proof of Proposition~\ref{prop:local-theory}.
\end{proof}

\paragraph{Continuous dependence on finite existence intervals.}\label{par:finite-interval-dependence}
Let $u$ be a solution in \eqref{eq:solution-class} on $[0,T]$ and
let $u_{0,n}\to u(0)$ in $H^s_\sigma$.
The compact set $\{u(t):0\le t\le T\}\subset H^s_\sigma$ is covered
by finitely many data neighborhoods from
Proposition~\ref{prop:local-theory}; the minimum of their common
positive lifespans is positive.  Subdivide $[0,T]$ into finitely
many intervals shorter than this minimum.  On each interval,
local continuous dependence and uniqueness identify the solution
starting from its left endpoint with the restriction of $u$.
Induction over the subdivision shows that, for all sufficiently
large $n$, the solution $u_n$ exists on $[0,T]$ and converges to $u$
in both local solution norms on every piece.  The global supremum
is the maximum of finitely many local suprema, and the square of
the dissipation norm is the sum of the corresponding squared
norms.  Hence the convergence holds in
$C([0,T];H^s)$ and in
$L^2(0,T;H^s)$ after applying $\Lambda_h^\alpha$.

\paragraph{Continuation with a small remainder.}\label{par:small-remainder-continuation}
The construction in the proof of Proposition~\ref{prop:local-theory}
also applies on any prescribed finite interval on which a fixed smooth
linear flow $v$ has bounded Fourier support and its remainder satisfies
the strict smallness condition used in the construction.
For completeness, suppose a solution on $[0,T_*)$, $T_*<\infty$,
has $u=v+w$, where $v=S_\Omega(t)v_0$ has bounded Fourier support and
\[
 \sup_{t<T_*}\|w(t)\|_{H^s}<c_{s,\alpha}\nu,\qquad
 \int_0^{T_*}\|\Lambda_h^\alpha w(t)\|_{H^s}^2\,dt<\infty.
\]
Choose $c_{s,\alpha}$ strictly below the smallness thresholds in
\eqref{eq:diff-common-background}--\eqref{eq:diff-alpha1} and
Lemma~\ref{lem:weighted-self}.
The weighted argument gives a uniform bound as in
\eqref{eq:construction-weighted-bound} up to $T_*$.
For a solution in the stated class, apply this argument first with
the bounded weights $\min(\omega_j,\omega_N)$ and then pass to the
limit by monotone convergence.
The equation bounds $\partial_tu$ in $H^{-m_0}$ for a sufficiently
large fixed integer $m_0$: the quadratic tensor is bounded in $L^1$,
and the linear multipliers have fixed finite order.
Hence $u(t)$ has a limit in this weaker space as $t\uparrow T_*$.
The weighted high-frequency tails and a fixed low-frequency cutoff
imply convergence in $H^s$.
Proposition~\ref{prop:local-theory}, applied to that limit, extends
the solution beyond $T_*$.  This is the continuation argument used
for the two global regimes below.
\subsection{Small-data global well-posedness}
\label{subsec:small-data-global}

Throughout the remainder of this section, $(\alpha,s)$ belongs to the
range of Theorem~\ref{thm:wellposed-main}.  The energy identity and
Proposition~\ref{prop:sobolev-self} give
\begin{equation}\label{eq:small-data-energy}
 \frac12\frac{d}{dt}\mathcal E_s(u)^2
 +\bigl(\nu-C_{s,\alpha}\mathcal E_s(u)\bigr)
                         \mathcal D_{s,\alpha}(u)^2\le0.
\end{equation}
Choose $c_{s,\alpha}>0$ so that
$\|u_0\|_{H^s}\le c_{s,\alpha}\nu$ places the initial energy strictly
below all the smallness thresholds in
Proposition~\ref{prop:local-theory}.
A continuity argument in \eqref{eq:small-data-energy} preserves this
smallness and gives
\[
 \sup_{t\ge0}\mathcal E_s(u(t))^2
       +\nu\int_0^\infty\mathcal D_{s,\alpha}(u(t))^2\,dt
       \le\mathcal E_s(u_0)^2.
\]
The construction in Proposition~\ref{prop:local-theory}, with zero
smooth background, applies on every finite time interval with this
same bound.  Uniqueness identifies the resulting solutions on
overlapping intervals.  Norm equivalence proves the small-data
assertion of Theorem~\ref{thm:wellposed-main}.  The threshold and the
energy bound are independent of $\Omega$, since the Coriolis term
has zero energy pairing in every scalar Fourier block.

\subsection{Frequency-localized estimates for the linear flow}
\label{subsec:rotation-linear}

Let
\[
 \mathcal G_\Omega(t)=e^{-t\Omega\mathbb P(e_3\times\cdot)}
\]
be the Coriolis group on divergence-free fields.  The semigroup
\eqref{eq:S-definition} factors as
$S_\Omega(t)=e^{-\nu t\Lambda_h^{2\alpha}}\mathcal G_\Omega(t)$.
The two factors commute, and $\mathcal G_\Omega(t)$ is unitary in
every Sobolev space.  To describe its phases, let
\[
 P(\xi)=I-\frac{\xi\otimes\xi}{|\xi|^2},
 \qquad
 \mathcal J(\xi)a=\frac{\xi}{|\xi|}\times a,
 \qquad
 P_\pm(\xi)=\tfrac12\bigl(P(\xi)\pm i\mathcal J(\xi)\bigr).
\]
On $\xi^\perp$,
$P(\xi)(e_3\times a)=(\xi_3/|\xi|)\mathcal J(\xi)a$ and
$\mathcal J(\xi)^2=-I$.  Consequently
\begin{equation}\label{eq:rotation-group-diagonal}
 \mathcal G_\Omega(t)
   =e^{i\Omega tD_3/|D|}\mathbb P_+
      +e^{-i\Omega tD_3/|D|}\mathbb P_-,
\end{equation}
where $\mathbb P_\pm$ have symbols $P_\pm$ and are bounded on $L^2$.

\begin{lemma}[Localized Coriolis estimate]
\label{lem:rotation-strichartz}
Let $\chi\in C_c^\infty(\mathbb R^3\setminus\{0\})$ and
$2<q<\infty$.  Then, for $\Omega\ne0$,
\begin{equation*}
 \|\chi(D)\mathcal G_\Omega(\cdot)f\|_{L^q(\mathbb R;L^\infty)}
       \le C_{\chi,q}|\Omega|^{-1/q}\|f\|_2
\end{equation*}
for every $f\in L^2_\sigma(\mathbb R^3)$.
\end{lemma}

\begin{proof}
The frequency-localized estimate of Koh, Lee and Takada
\cite{KohLeeTakada2014a}, in the form stated in
\cite[Theorem~4.1]{Takada2014Survey}, is
\[
 \|\psi(D)e^{\pm itD_3/|D|}f\|_{L^q(\mathbb R;L^p)}
       \le C_{\psi,q,p}\|f\|_2
\]
for a fixed annular cutoff $\psi$, provided
\[
 2\le q,p\le\infty,\qquad
 \frac1q+\frac1p\le\frac12,\qquad(q,p)\ne(2,\infty).
\]
For a dyadic annulus of radius $2^j$, the degree-zero homogeneity of
the phase and spatial scaling multiply the right side by
$2^{3j(1/2-1/p)}$.  A finite annular partition therefore gives the
same estimate for $\chi$, with a constant depending on its fixed
support.  Take $p=\infty$ and $q>2$, then change variables from $t$
to $\Omega t$.  Applying the resulting scalar bound to
$\mathbb P_+f$ and $\mathbb P_-f$ in
\eqref{eq:rotation-group-diagonal} proves the lemma.
The factor $|\Omega|^{-1/q}$ comes from this time rescaling.  Since
$p=\infty$ requires $q>2$, it gives precisely the range
$\gamma=1/q<1/2$ in
Theorem~\ref{thm:wellposed-main}\textup{(iv)}.
\end{proof}

\begin{lemma}[Spectral approximation]
\label{lem:rotation-spectral-approximation}
Let $s\ge0$ and $u_0\in H^s_\sigma(\mathbb R^3)$ be real-valued.
For every $\eta>0$ there is a real even
$\chi\in C_c^\infty(\mathbb R^3)$, with $0\le\chi\le1$, such that
\begin{equation*}
 \begin{gathered}
 v_0=\chi(D)u_0,\qquad w_0=(1-\chi(D))u_0,\qquad
                  \|w_0\|_{H^s}<\eta,\\
 \operatorname{supp}\widehat v_0
     \subset\{\xi:|\xi|\le R,\ |\xi_h|\ge\kappa\}
 \end{gathered}
\end{equation*}
for some $R>\kappa>0$.  Both fields are real and divergence-free.
\end{lemma}

\begin{proof}
Choose a smooth nondecreasing function $\zeta$ on $[0,\infty)$,
equal to zero on $[0,1]$ and to one on $[2,\infty)$, and put
\[
 \chi_n(\xi)=\zeta(n|\xi_h|)
                     \bigl(1-\zeta(|\xi|/n)\bigr).
\]
These are real even smooth cutoffs supported where
$|\xi_h|\ge n^{-1}$ and $|\xi|\le2n$.
They tend to one almost everywhere, since the vertical frequency
axis has measure zero.  Dominated convergence in the Fourier
definition of $H^s$ gives
$\|(1-\chi_n(D))u_0\|_{H^s}\to0$.
Scalar Fourier multipliers preserve divergence freeness, and
reality and evenness of $\chi_n$ preserve real-valued fields.
\end{proof}

\begin{proposition}[Integrable bounds for the localized linear flow]
\label{prop:rotation-localized-linear}
Suppose
\[
 \operatorname{supp}\widehat v_0
       \subset\{\xi:|\xi|\le R,\ |\xi_h|\ge\kappa\},
 \qquad R>\kappa>0,
\]
and let $v=S_\Omega(\cdot)v_0$.  For every integer $k\ge0$ and
$2<q<\infty$,
\begin{align}
 \|v(t)\|_{H^k}
   &\le C_{k,R}e^{-\nu\kappa^{2\alpha}t}\|v_0\|_2,
       \label{eq:rotation-linear-decay}\\
 \|v\|_{L^q(0,\infty;W^{k,\infty})}
   +\|v\|_{L^1(0,\infty;W^{k,\infty})}
   &\le C_{k,q,R,\kappa,\nu,\alpha}
                          |\Omega|^{-1/q}\|v_0\|_2,
       \qquad \Omega\ne0.
       \label{eq:rotation-linear-integrable}
\end{align}
For every $s>0$ and integer $k>s+1$,
\begin{equation}\label{eq:rotation-linear-forcing}
 \|v\cdot\nabla v\|_{L^1(0,\infty;H^s)}
       \le C_{s,k}\|v_0\|_{H^s}
                   \|v\|_{L^1(0,\infty;W^{k,\infty})}.
\end{equation}
\end{proposition}

\begin{proof}
The Fourier support is preserved by the flow.
Unitarity of $\mathcal G_\Omega(t)$ and the lower bound on
$|\xi_h|$ give
$\|v(t)\|_2\le e^{-\nu\kappa^{2\alpha}t}\|v_0\|_2$.
Bernstein's inequality yields \eqref{eq:rotation-linear-decay}.

For $0<\alpha\le1$, the horizontal fractional heat semigroup is
contractive on $L^\infty$.  Choose
$\widetilde\chi\in C_c^\infty(\mathbb R^3\setminus\{0\})$ equal
to one on the specified spectral support.  For every
$|\beta|\le k$,
\[
 \|\partial^\beta v(t)\|_\infty
 \le\|\widetilde\chi(D)\mathcal G_\Omega(t)
                                  \partial^\beta v_0\|_\infty.
\]
Lemma~\ref{lem:rotation-strichartz} and
$\|\partial^\beta v_0\|_2\le C_{k,R}\|v_0\|_2$ give the
$L^q$ bound in \eqref{eq:rotation-linear-integrable}.
For the $L^1$ bound, use the semigroup property on each unit
interval:
\[
 \begin{aligned}
 \int_n^{n+1}\|v(t)\|_{W^{k,\infty}}\,dt
 &\le
 \|S_\Omega(\cdot)S_\Omega(n)v_0
                         \|_{L^q(0,1;W^{k,\infty})}\\
 &\le C|\Omega|^{-1/q}
             e^{-\nu\kappa^{2\alpha}n}\|v_0\|_2 .
 \end{aligned}
\]
Summing over $n\ge0$ proves the assertion.

Multiplication by a function in $W^{k-1,\infty}$ is bounded on
$H^s$ when $k-1>s$.  Applying this to the factors
$\partial_jv$ gives
\[
 \|v\cdot\nabla v\|_{H^s}
     \le C_{s,k}\|v\|_{H^s}\|v\|_{W^{k,\infty}}.
\]
Since $S_\Omega(t)$ is contractive on $H^s$, integration proves
\eqref{eq:rotation-linear-forcing}.
\end{proof}

\subsection{A quantitative estimate for the nonlinear remainder}
\label{subsec:rotation-remainder}

Fix an integer $k=k(s)>s+3$.  This stronger choice than in
Proposition~\ref{prop:rotation-localized-linear} allows the transport
commutator and the high-frequency coefficient terms to be estimated
with the same norm.  Let $v=S_\Omega(\cdot)v_0$ be a smooth linear flow
and write $u=v+w$.  The remainder solves
\begin{equation*}
 \begin{aligned}
 \partial_tw+\nu\Lambda_h^{2\alpha}w
 &+\Omega\mathbb P(e_3\times w)\\
 &+\mathbb P(w\cdot\nabla w+v\cdot\nabla w+w\cdot\nabla v)
       =-\mathbb P(v\cdot\nabla v).
 \end{aligned}
\end{equation*}
The smooth-coefficient terms satisfy
\[
 \left|\sum_{j\in\mathbb Z}2^{2js}
   \bigl(\Delta_j(v\cdot\nabla w+w\cdot\nabla v),\Delta_jw\bigr)_{L^2}
       +(v\cdot\nabla w+w\cdot\nabla v,w)_{L^2}\right|
       \le C_{s,k}\|v\|_{W^{k,\infty}}\mathcal E_s(w)^2.
\]
Indeed, the low-frequency transport cancels by
$\operatorname{div}v=0$, and its commutator is controlled by
the first moment of the dyadic kernel and $\|\nabla v\|_\infty$.
For a high-frequency coefficient, the bound
$\|\Delta_jv\|_\infty\le C2^{-kj}\|v\|_{W^{k,\infty}}$ for
$j\ge0$ makes the remaining sums convergent.  Negative indices are
controlled by the $L^2$ part of $\mathcal E_s$ and $s>0$.
The product $w\cdot\nabla v$ is bounded by the $H^s$ multiplier
estimate used to prove \eqref{eq:rotation-linear-forcing}.
The same argument applies to the frequency weights satisfying the adjacent-ratio bound in
Lemma~\ref{lem:weighted-self}.
Proposition~\ref{prop:sobolev-self} consequently gives
\begin{equation}\label{eq:rotation-remainder-energy}
 \begin{aligned}
 \frac12\frac{d}{dt}\mathcal E_s(w)^2
  +\bigl(\nu-C_{s,\alpha}\mathcal E_s(w)\bigr)
                              \mathcal D_{s,\alpha}(w)^2
 \le{}& C_{s,k}\|v\|_{W^{k,\infty}}\mathcal E_s(w)^2\\
 &+C_s\|v\cdot\nabla v\|_{H^s}\mathcal E_s(w).
 \end{aligned}
\end{equation}

\begin{lemma}[Global remainder estimate]
\label{lem:rotation-global-remainder}
There are constants $c_{s,\alpha}>0$ and $C>0$, depending only
on $s,\alpha$ and the fixed integer $k$, with the following property.
Let $v_0\in L^2_\sigma$ have bounded Fourier support,
$w_0\in H^s_\sigma$, and $v=S_\Omega(\cdot)v_0$.
Suppose
\begin{equation}\label{eq:rotation-remainder-smallness}
 \left(\mathcal E_s(w_0)
       +C\|v\cdot\nabla v\|_{L^1(0,\infty;H^s)}\right)
 \exp\!\left(C\|v\|_{L^1(0,\infty;W^{k,\infty})}\right)
                   \le c_{s,\alpha}\nu.
\end{equation}
Then the solution with initial datum $v_0+w_0$ is global, and
\begin{equation}\label{eq:rotation-remainder-sup}
 \sup_{t\ge0}\mathcal E_s(w(t))
 \le
 \left(\mathcal E_s(w_0)
       +C\|v\cdot\nabla v\|_{L^1(0,\infty;H^s)}\right)
 \exp\!\left(C\|v\|_{L^1(0,\infty;W^{k,\infty})}\right).
\end{equation}
Moreover,
\begin{equation}\label{eq:rotation-remainder-dissipation}
 \begin{aligned}
 \nu\int_0^\infty\mathcal D_{s,\alpha}(w)^2\,dt
 \le{}&\mathcal E_s(w_0)^2
  +2C\|v\|_{L^1(0,\infty;W^{k,\infty})}
                      \sup_{t\ge0}\mathcal E_s(w(t))^2\\
 &+2C\|v\cdot\nabla v\|_{L^1(0,\infty;H^s)}
                      \sup_{t\ge0}\mathcal E_s(w(t)).
 \end{aligned}
\end{equation}
If also $\|v\|_{L^1(0,\infty;W^{k,\infty})}\le1$, then
\begin{equation}\label{eq:rotation-remainder-full-bound}
 \begin{aligned}
 &\sup_{t\ge0}\|w(t)\|_{H^s}
   +\sqrt\nu\,
       \|\Lambda_h^\alpha w\|_{L^2(0,\infty;H^s)}\\
 &\hspace{16mm}\le
 C\left(\|w_0\|_{H^s}
             +\|v\cdot\nabla v\|_{L^1(0,\infty;H^s)}\right).
 \end{aligned}
\end{equation}
All constants in these estimates are independent of the spectral
cutoff, the datum, and $\Omega$.
\end{lemma}

\begin{proof}
Decrease $c_{s,\alpha}$ to include the smallness thresholds in
Proposition~\ref{prop:local-theory}, Lemma~\ref{lem:weighted-self},
and the higher-order bound \eqref{eq:higher-sobolev-self}; for
$\alpha=1$, use \eqref{eq:laplacian-self-mixed} and Young's inequality.
On any interval on which
$\mathcal E_s(w)\le2c_{s,\alpha}\nu$, the coefficient of
$\mathcal D_{s,\alpha}(w)^2$ in
\eqref{eq:rotation-remainder-energy} is at least $\nu/2$.
Divide the resulting scalar inequality by the energy norm after
replacing it by
$\bigl(\mathcal E_s(w)^2+\epsilon^2\bigr)^{1/2}$, and then let
$\epsilon\downarrow0$.  Gronwall's inequality gives
\begin{equation}\label{eq:remainder-bootstrap}
 \begin{aligned}
 \mathcal E_s(w(t))
 \le{}&
 \left(\mathcal E_s(w_0)
       +C\int_0^t\|v\cdot\nabla v\|_{H^s}\,d\tau\right)\\
 &\qquad\times
       \exp\!\left(C\int_0^t\|v\|_{W^{k,\infty}}\,d\tau\right)
 \le c_{s,\alpha}\nu.
 \end{aligned}
\end{equation}
This gives a strict improvement of the bootstrap bound.
Integrating the absorbed energy inequality gives
\begin{equation}\label{eq:remainder-finite-energy}
 \begin{aligned}
 \mathcal E_s(w(T))^2
       +\nu\int_0^T\mathcal D_{s,\alpha}(w)^2\,dt
 \le{}&\mathcal E_s(w_0)^2
       +2C\int_0^T\|v\|_{W^{k,\infty}}\mathcal E_s(w)^2\,dt\\
 &+2C\int_0^T\|v\cdot\nabla v\|_{H^s}\mathcal E_s(w)\,dt .
 \end{aligned}
\end{equation}
On the current existence interval,
\eqref{eq:remainder-bootstrap}--\eqref{eq:remainder-finite-energy}
give a strict improvement of the remainder bound and finite total dissipation
on every finite time interval.  The continuation argument following
Proposition~\ref{prop:local-theory}
(p.~\pageref{par:small-remainder-continuation}) therefore excludes
a finite maximal existence time.  The weighted dyadic estimates give a
strong $H^s$ trace and convergence of the approximations.
Uniqueness identifies the solutions constructed on different finite
intervals.  Hence the solution is global and
\eqref{eq:rotation-remainder-sup} holds.
Taking suprema on the right and letting $T\to\infty$ proves
\eqref{eq:rotation-remainder-dissipation}.

Under the additional unit bound on
$\|v\|_{L^1W^{k,\infty}}$, \eqref{eq:rotation-remainder-sup} gives
\[
 \sup_{t\ge0}\mathcal E_s(w(t))
 \le C\left(\mathcal E_s(w_0)
                 +\|v\cdot\nabla v\|_{L^1H^s}\right).
\]
The last term of \eqref{eq:rotation-remainder-dissipation} is
therefore bounded by a constant times
$\bigl(\mathcal E_s(w_0)+\|v\cdot\nabla v\|_{L^1H^s}\bigr)^2$.
The same is true of the first two terms.  Taking square roots and
using norm equivalence proves
\eqref{eq:rotation-remainder-full-bound}.
\end{proof}

\subsection{Rapid rotation and approximation by the linear flow}
\label{subsec:rotation-global-proof}

Fix $\nu,\alpha,s$ and one datum $u_0\in H^s_\sigma$.
Choose $\eta>0$ small enough that
$\mathcal E_s(w_0)\le C_s\eta<c_{s,\alpha}\nu/2$
in Lemma~\ref{lem:rotation-spectral-approximation}, and keep the
resulting cutoff fixed.  Proposition~\ref{prop:rotation-localized-linear}
gives, for $v^\Omega=S_\Omega(\cdot)v_0$,
\begin{equation}\label{eq:rotation-fixed-cutoff-smallness}
 \|v^\Omega\|_{L^1(0,\infty;W^{k,\infty})}
 +\|v^\Omega\cdot\nabla v^\Omega\|_{L^1(0,\infty;H^s)}
       \longrightarrow0
       \qquad (|\Omega|\to\infty).
\end{equation}
Thus \eqref{eq:rotation-remainder-smallness} holds for all
sufficiently large $|\Omega|$, and
Lemma~\ref{lem:rotation-global-remainder} gives the global solution.
The threshold can depend on the chosen cutoff and hence on the
frequency distribution of the initial datum.
The linear energy identity
\begin{equation}\label{eq:rotation-linear-energy}
 \|S_\Omega(T)g\|_{H^s}^2
 +2\nu\int_0^T
       \|\Lambda_h^\alpha S_\Omega(t)g\|_{H^s}^2\,dt
       =\|g\|_{H^s}^2
\end{equation}
and \eqref{eq:rotation-remainder-full-bound} show that the global
solution has bounded $H^s$ norm and finite total dissipation.

To prove the linear-approximation assertion, take an arbitrarily
small $\eta$ as above and write the resulting decomposition as
$u_0=v_{0,\eta}+w_{0,\eta}$.
Uniqueness identifies the global solution constructed from this
decomposition with $u^\Omega$ whenever $|\Omega|$ is sufficiently
large.  Apply \eqref{eq:rotation-remainder-full-bound} to
$u^\Omega-S_\Omega(\cdot)v_{0,\eta}$ and
\eqref{eq:rotation-linear-energy} to $w_{0,\eta}$.  The triangle
inequality gives
\begin{equation*}
 \begin{aligned}
 &\sup_{t\ge0}
     \|u^\Omega(t)-S_\Omega(t)u_0\|_{H^s}\\
 &\quad+\sqrt\nu\,
   \|\Lambda_h^\alpha(u^\Omega-S_\Omega(\cdot)u_0)
                              \|_{L^2(0,\infty;H^s)}\\
 &\qquad\le
 C\eta+
 C\|S_\Omega(\cdot)v_{0,\eta}\cdot
            \nabla S_\Omega(\cdot)v_{0,\eta}
                              \|_{L^1(0,\infty;H^s)}.
 \end{aligned}
\end{equation*}
The constant $C$ is independent of $\eta$ and the cutoff, whereas
the required rotation threshold may depend on both.
First let $|\Omega|\to\infty$ at fixed $\eta$, using
\eqref{eq:rotation-fixed-cutoff-smallness}.  Then let
$\eta\downarrow0$.  This proves convergence to the full linear
Stokes--Coriolis evolution in both global norms in
Theorem~\ref{thm:wellposed-main}.

\subsection{Rate for spectrally localized data}
\label{subsec:rotation-rate}

Assume now that
\[
 \operatorname{supp}\widehat u_0
       \subset\{\xi:|\xi|\le R,\ |\xi_h|\ge\kappa\}.
\]
Take $v_0=u_0$, $w_0=0$ and
$v^\Omega=S_\Omega(\cdot)u_0$.
Given $0<\gamma<1/2$, choose $q=1/\gamma>2$ in
Proposition~\ref{prop:rotation-localized-linear}.  Then
\[
 \begin{aligned}
 \|v^\Omega\|_{L^1(0,\infty;W^{k,\infty})}
       &\le C_\gamma|\Omega|^{-\gamma}\|u_0\|_2,\\
 \|v^\Omega\cdot\nabla v^\Omega\|_{L^1(0,\infty;H^s)}
       &\le C_\gamma|\Omega|^{-\gamma}
                                \|u_0\|_{H^s}\|u_0\|_2.
 \end{aligned}
\]
For all sufficiently large $|\Omega|$, the hypotheses of
Lemma~\ref{lem:rotation-global-remainder} hold, including the unit
bound $\|v^\Omega\|_{L^1(0,\infty;W^{k,\infty})}\le1$.  Estimate
\eqref{eq:rotation-remainder-full-bound}, with $w_0=0$, gives
\[
 \begin{aligned}
 &\sup_{t\ge0}
       \|u^\Omega(t)-S_\Omega(t)u_0\|_{H^s}\\
 &\qquad+\sqrt\nu\,
  \|\Lambda_h^\alpha(u^\Omega-S_\Omega(\cdot)u_0)
                                  \|_{L^2(0,\infty;H^s)}
       \le C_\gamma|\Omega|^{-\gamma}.
 \end{aligned}
\]
The constants and the threshold depend only on
$\gamma,s,\alpha,\nu,R,\kappa$ and $\|u_0\|_{H^s}$.
This proves the final assertion of
Theorem~\ref{thm:wellposed-main}.

\section{Norm inflation below the critical index}
\label{sec:subcritical-inflation}

We first prove Theorem~\ref{thm:illposed-main} for
$0<s<5/2-2\alpha$.  The argument transfers a fixed smooth Euler
solution to \eqref{eq:NSC} through a scaling under which both the
horizontal viscosity and the rotation parameter tend to zero.

\begin{proof}[Proof of Theorem~\ref{thm:illposed-main} for
$0<s<5/2-2\alpha$]
Fix $\nu>0$, $\Omega\in\mathbb R$, and $\varepsilon,B,T_0>0$.
By Luo's Euler norm-inflation theorem
\cite[Theorem~1.1]{Luo2026Sobolev}, there are real divergence-free
$U_0\in C_c^\infty(\mathbb R^3)$ and a smooth Euler solution $U$ on
a closed interval $[0,T]$ such that
\begin{equation*}
 \|U_0\|_{H^s}<\varepsilon,\qquad
 \|U(T)\|_{\dot H^s}>B+2.
\end{equation*}
The small parameter in that theorem can be chosen to meet these two
inequalities simultaneously.  All Sobolev norms of this fixed
solution are bounded on $[0,T]$.

For $N\ge1$, set
\begin{equation*}
 \nu_N=\nu N^{s+2\alpha-5/2},\qquad
 \Omega_N=\Omega N^{s-5/2}.
\end{equation*}
Both coefficients tend to zero.  Let $V_N$ be the classical local
solution of
\[
 \partial_\tau V_N+\nu_N\Lambda_h^{2\alpha}V_N
 +\Omega_N\mathbb P(e_3\times V_N)
 +\mathbb P(V_N\cdot\nabla V_N)=0,\qquad V_N(0)=U_0.
\]
We claim that, for sufficiently large $N$, this solution exists
through $T$ and converges to $U$ in $C([0,T];H^3)$.

Indeed, $R_N=V_N-U$ satisfies
\begin{align*}
 \partial_\tau R_N+\nu_N\Lambda_h^{2\alpha}R_N
 +\Omega_N\mathbb P(e_3\times R_N)
 &=-\mathbb P\bigl((U+R_N)\cdot\nabla R_N
                                  +R_N\cdot\nabla U\bigr)\\
 &\quad-\nu_N\Lambda_h^{2\alpha}U
   -\Omega_N\mathbb P(e_3\times U).
\end{align*}
The integer Sobolev commutator estimate and
$H^3\hookrightarrow W^{1,\infty}$ give
\begin{align*}
 \frac12\frac{d}{d\tau}\|R_N\|_{H^3}^2
 +\nu_N\|\Lambda_h^\alpha R_N\|_{H^3}^2
 &\le C\bigl(\|U\|_{H^4}+\|R_N\|_{H^3}\bigr)
                                      \|R_N\|_{H^3}^2\\
 &\quad+
 \bigl(\nu_N\|U\|_{H^5}+|\Omega_N|\|U\|_{H^3}\bigr)
                                      \|R_N\|_{H^3}.
\end{align*}
The Coriolis term paired with $R_N$ vanishes.  We also used
$\|\Lambda_h^{2\alpha}U\|_{H^3}\le\|U\|_{H^5}$, since $\alpha\le1$.
On the interval where $\|R_N\|_{H^3}\le1$, Gronwall's inequality
therefore yields
\begin{equation}\label{eq:negative-error}
 \sup_{0\le\tau\le T}\|R_N(\tau)\|_{H^3}
          \le C_{U,T}\bigl(\nu_N+|\Omega_N|\bigr).
\end{equation}
Initially the supremum is restricted to the common existence
interval.  For large $N$, the right side is smaller than $1/2$.
The classical $H^3$ continuation criterion then extends $V_N$
through $T$, and higher regularity follows from the corresponding
energy estimates.  This proves the claim.

Define
\begin{equation*}
 u_N(t,x)=N^{3/2-s}V_N(N^{5/2-s}t,Nx).
\end{equation*}
The homogeneity of the Leray projector and of
$\Lambda_h^{2\alpha}$ shows that $u_N$ solves \eqref{eq:NSC} with
the original coefficients $\nu$ and $\Omega$.  Its initial datum
is real, divergence-free, and compactly supported.  Fourier change
of variables gives
\begin{align*}
 \|u_N(0)\|_{H^s}^2
 &=\int_{\mathbb R^3}(N^{-2}+|\xi|^2)^s
                       |\widehat U_0(\xi)|^2\,d\xi
 \le\|U_0\|_{H^s}^2,\\
 \|N^{3/2-s}F(N\cdot)\|_{\dot H^s}
 &=\|F\|_{\dot H^s}.
\end{align*}
The physical times $t_N=TN^{s-5/2}$ tend to zero.  Since
$s<5/2<3$, \eqref{eq:negative-error} implies
\[
 \|u_N(t_N)\|_{\dot H^s}
 =\|V_N(T)\|_{\dot H^s}
 \ge\|U(T)\|_{\dot H^s}-\|V_N(T)-U(T)\|_{H^3}>B
\]
for all sufficiently large $N$.  Taking also $t_N<T_0$ proves
the assertion in this range.
\end{proof}

\section{Norm inflation at the horizontal-Laplacian endpoint}
\label{sec:endpoint-inflation}

It remains to prove Theorem~\ref{thm:illposed-main} when
$\alpha=1$ and $s=1/2$.  We first construct a finite-energy
approximate solution at zero rotation and control its error uniformly
in the vertical concentration scales.  A critical scaling then
transfers the construction to every prescribed finite $\Omega$.
Throughout this section, $z=x_3$ and $\nu>0$ is fixed.

\subsection{Conormal analytic estimates}
\label{subsec:endpoint-calculus}

For $m\ge0$, define
\[
 \|f\|_{\mathcal A_m}
 =\int_{\mathbb R^2}(1+|\xi_h|)^m
          \|\mathcal F_hf(\xi_h,\cdot)\|_{L_z^\infty}\,d\xi_h,
 \qquad Z=z\partial_z,
\]
Here $Z$ is the vertical dilation vector field.  For functions of
$x_h$ alone, $\mathcal A_m$ denotes the same integral without the
vertical supremum.  Set
\[
 \|f\|_{X_{\rho,m}}
 =\sum_{n\ge0}\frac{\rho^n}{n!}\|Z^nf\|_{\mathcal A_m},
 \qquad
 Y_{\rho,m}(f)=\|Zf\|_{X_{\rho,m}}
                    =\partial_\rho\|f\|_{X_{\rho,m}}.
\]
Component norms are summed for vector fields.  We suppress $m$
when it is fixed.  The order of the horizontal integral and the
vertical supremum will be used in the pressure estimates.

\begin{lemma}[Analytic operator estimates]
\label{lem:endpoint-analytic-calculus}
There is $\rho_*>0$ such that, for $0\le\rho\le\rho_*$ and
$0<\eta\le1$, the operator $J_\eta=e^{\eta\Delta}$ satisfies
\begin{align}
 \|fg\|_{X_\rho}
   &\le C_m\|f\|_{X_\rho}\|g\|_{X_\rho},
 \notag\\
 \|\mathbb Pf\|_{X_\rho}+\|J_\eta f\|_{X_\rho}
   &\le C\|f\|_{X_\rho},
 \label{eq:analytic-stability-operators}\\
 Y_\rho(J_\eta f)
   &\le C\bigl(Y_\rho(f)+\|f\|_{X_\rho}\bigr),
 \label{eq:analytic-stability-J-Y}\\
 \|\nabla J_\eta f\|_{X_\rho}
   &\le C\eta^{-1/2}\|f\|_{X_\rho}.
 \label{eq:Gaussian-gradient-bound}
\end{align}
Moreover,
\begin{equation}\label{eq:analytic-stability-interpolation}
 \|\Lambda_hf\|_{X_\rho}^2
 \le\|f\|_{X_\rho}\|\Lambda_h^2f\|_{X_\rho}.
\end{equation}
If $f$ is divergence-free, $f_h$ is even in $z$, and $f_3$ is odd,
then
\begin{equation}\label{eq:analytic-stability-c}
 c_f:=\frac{f_3}{z}
 =-\int_0^1\operatorname{div}_h f_h(x_h,az)\,da,
 \qquad \|c_f\|_{X_\rho}\le C\|\Lambda_hf\|_{X_\rho},
\end{equation}
with continuous extension at zero, and
\begin{equation}\label{eq:analytic-stability-c-J}
             \|c_{J_\eta f}\|_{X_\rho}\le C\|c_f\|_{X_\rho}.
\end{equation}
All constants are independent of $\eta$.
\end{lemma}

\begin{proof}
The product estimate follows from horizontal convolution, the
submultiplicativity of $(1+|\xi_h|)^m$, and the factorial Leibniz
formula.  Inequality \eqref{eq:analytic-stability-interpolation}
is Cauchy--Schwarz for the positive measure consisting of the
horizontal frequency integral and the sum over conormal orders.

At fixed $k=|\xi_h|>0$, the kernel of
$(k^2-\partial_z^2)^{-1}$ is
\[
                         \frac1{2k}e^{-k|z|}.
\]
The entries of $\mathbb P$ are combinations of the identity and
convolution by
\[
 k e^{-k|z|},\qquad k\operatorname{sgn}(z)e^{-k|z|},
\]
with bounded coefficients depending on $\xi_h/k$.  Their $L_z^1$
norms are uniform in $k$.  For any convolution kernel $K$,
\[
                         [Z,K*]=(K+zK')*.
\]
For either exponential kernel, the real dilation family
$e^t kH(ke^t z)$ extends to an $L^1$-valued holomorphic function
for $|\operatorname{Im}t|<\pi/2$.  For the odd kernel this
extension is
$e^t k\operatorname{sgn}(z)e^{-ke^t|z|}$.
Its $L^1$ norm is uniformly bounded on each fixed disk strictly
inside that strip.  Cauchy's estimate gives
\[
 \sup_{k>0}\|(1+z\partial_z)^r[kH(kz)]\|_1
                         \le C_0C_1^r r!.
\]
The identity kernel has zero commutator.  Expanding $Z^n\mathbb P$
by commutators and summing gives
\[
 \|\mathbb Pf\|_{X_\rho}
 \le\frac{C_0}{1-C_1\rho}\|f\|_{X_\rho}.
\]
The horizontal frequency $k=0$ is a null set in this norm and can
also be assigned its vertical Fourier multiplier value.

After horizontal Fourier transformation, $J_\eta$ is multiplication
by $e^{-\eta|\xi_h|^2}$ followed by convolution with
\[
                  g_\eta(z)=(4\pi\eta)^{-1/2}e^{-z^2/(4\eta)}.
\]
The analytic dilation family $e^tg_\eta(e^tz)$ has uniformly bounded
$L^1$ norm on each fixed disk $|t|<a<\pi/4$.  It therefore gives
the same factorial commutator estimate, uniformly in $\eta$, and
proves the Gaussian part of \eqref{eq:analytic-stability-operators}.
The commutators of $[Z,J_\eta]$ have the bounds
$C_0C_1^{r+1}(r+1)!$; their weighted series is convergent for small
$\rho$.  This proves \eqref{eq:analytic-stability-J-Y}.
For a vertical derivative, the same argument applied to
$g_\eta'$ gives the additional factor $\eta^{-1/2}$.
For horizontal derivatives use
$|\xi_h|e^{-\eta|\xi_h|^2}\le C\eta^{-1/2}$.
This proves \eqref{eq:Gaussian-gradient-bound}.

Formula \eqref{eq:analytic-stability-c} follows by integrating
$\partial_zf_3=-\operatorname{div}_hf_h$ from zero.
The resulting dilation integral commutes with $Z$ and is bounded
on vertical $L^\infty$.

To prove \eqref{eq:analytic-stability-c-J}, let $c$ be even and set
\[
                  T_\eta c(z)=\frac{e^{\eta\partial_z^2}(zc)(z)}{z}.
\]
For $z>0$, the heat kernel gives
\[
 T_\eta c(z)
 =\frac1{z\sqrt{4\pi\eta}}\int_0^\infty r
 \left(e^{-(z-r)^2/(4\eta)}-e^{-(z+r)^2/(4\eta)}\right)c(r)\,dr.
\]
Angular integration identifies this with the three-dimensional heat
semigroup acting on the radial function $y\mapsto c(|y|)$.
Conjugating this semigroup by dilation gives the kernel
\[
 H_{\eta,t}(y)
 =e^{3t}(4\pi\eta)^{-3/2}
                 \exp\!\left(-\frac{e^{2t}|y|^2}{4\eta}\right),
 \qquad
 \|H_{\eta,t}\|_1
       =\bigl(\cos(2\operatorname{Im}t)\bigr)^{-3/2}.
\]
For a fixed $a<\pi/4$, Cauchy's estimate bounds the
$r$-fold commutator of $T_\eta$ with $Z$ in
$L^\infty\to L^\infty$ operator norm by $C_a r!a^{-r}$.
The factorial sum then implies
$\|T_\eta c\|_{X_\rho}\le C\|c\|_{X_\rho}$.
Finally,
\[
                  c_{J_\eta f}=e^{\eta\Delta_h}T_\eta c_f,
\]
and horizontal heat is contractive in $X_\rho$.
\end{proof}

\begin{lemma}[Strong convergence of Gaussian regularization]
\label{lem:gaussian-conormal-regularizer}
Let $0\le\rho<\rho_1\le\rho_*$.
If $f$ is a Schwartz function with finite $X_{\rho_1,m}$ norm, then
\[
                         \|J_\eta f-f\|_{X_{\rho,m}}\longrightarrow0
                         \qquad(\eta\to0).
\]
The convergence is uniform for families continuous in the Schwartz
topology on a compact parameter set with a uniform
$X_{\rho_1,m}$ bound.  After a further arbitrarily small decrease of the analytic radius,
the convergence also holds after one additional $Z$ derivative.
\end{lemma}

\begin{proof}
For finitely many conormal derivatives, the Gaussian kernel is an
approximate identity.  Its positive-order commutator kernels have
zero integral, uniformly bounded $L^1$ norms, and mass concentrating
at zero; they converge to zero on each fixed uniformly continuous
profile.  Dominated convergence handles the horizontal integral.
The remaining conormal sum is bounded by a constant times
$(\rho/\rho_1)^{n+1}\|f\|_{X_{\rho_1,m}}$, using
Lemma~\ref{lem:endpoint-analytic-calculus} at $\rho_1$.
This proves strong convergence.  The argument is uniform over the
stated compact families.  The additional derivative follows from
\[
 \|Zg\|_{X_{\rho,m}}
             \le\frac1{\rho'-\rho}\|g\|_{X_{\rho',m}},
 \qquad \rho<\rho'<\rho_1.
\]
\end{proof}

\subsection{Stability of an approximate solution}
\label{sec:analytic-conormal-stability}

\begin{lemma}[Smooth Gaussian approximations]
\label{lem:friedrichs-global-smooth}
For every fixed $\eta>0$ and
$u_0\in\bigcap_{k\ge0}H^k_\sigma(\mathbb R^3)$, the equation
\[
 \partial_tu^\eta-\nu\Delta_hu^\eta
 +J_\eta\mathbb P\bigl((J_\eta u^\eta)\cdot
                              \nabla(J_\eta u^\eta)\bigr)=0,
 \qquad u^\eta(0)=u_0,
\]
has a global smooth solution and satisfies the $L^2$ energy
identity.  It preserves the reflection parity in
\eqref{eq:analytic-stability-c}.  If $u_0\in X_{\rho,m}$, it also
has a local mild solution in $X_{\rho,m}$ agreeing with the
smooth solution.
\end{lemma}

\begin{proof}
For fixed $\eta$, the regularized nonlinearity is locally Lipschitz
on $L^2$.  Put $w=J_\eta u^\eta$.  Self-adjointness, commutation
with the Leray projector, and $\operatorname{div}w=0$ give
\[
 \left\langle J_\eta\mathbb P(w\cdot\nabla w),u^\eta\right\rangle
                       =\langle w\cdot\nabla w,w\rangle=0.
\]
Thus
\[
 \frac12\frac{d}{dt}\|u^\eta\|_2^2
                +\nu\|\nabla_hu^\eta\|_2^2=0,
\]
which extends the $L^2$ solution globally.  For every integer
$k\ge3$, ordinary derivatives commute with $J_\eta$, so
\[
 \frac{d}{dt}\|u^\eta\|_{H^k}^2
 +2\nu\|\nabla_hu^\eta\|_{H^k}^2
 \le C_k\|\nabla J_\eta u^\eta\|_\infty\|u^\eta\|_{H^k}^2.
\]
The bound
$\|\nabla J_\eta u^\eta\|_\infty
\le C\eta^{-5/4}\|u_0\|_2$ proves global smoothness.
Parity follows from uniqueness and invariance of the equation.
Finally, Lemma~\ref{lem:endpoint-analytic-calculus} gives
\[
 \|J_\eta\mathbb P((J_\eta u)\cdot\nabla J_\eta u)\|_{X_\rho}
                  \le C\eta^{-1/2}\|u\|_{X_\rho}^2
\]
and the corresponding local Lipschitz estimate.
Duhamel's formula for the horizontal heat semigroup gives a local mild solution in the
fixed space $X_\rho$.  Constructing it in $X_\rho\cap L^2$
identifies it with the smooth solution.
\end{proof}

\begin{proposition}[Stability in a conormal analytic norm]
\label{prop:endpoint-analytic-stability}
Fix $T>0$, $\nu>0$, $m\geq0$, and a sufficiently small $\rho_0>0$.
Let $U$ be a real smooth divergence-free field on
$[0,T]\times\mathbb R^3$ with $U_h$ even and $U_3$ odd in $z$.
Assume $U$ and $\partial_tU$ depend continuously on $t$ in the
Schwartz topology, and that
\begin{equation}\label{eq:analytic-stability-A}
 \sup_{0\leq t\leq T}
 \left(\|U(t)\|_{X_{2\rho_0,m}}
       +\|\Lambda_hU(t)\|_{X_{2\rho_0,m}}
       +Y_{2\rho_0,m}(U(t))\right)\leq A.
\end{equation}
The radius $2\rho_0$ is taken below the operator radius $\rho_*$.
Assume also that
\begin{equation}\label{eq:analytic-stability-v0}
                         \|\partial_zU(0)\|_{X_{\rho_0,m}}<\infty.
\end{equation}
Define the unprojected residual
\[
 \mathcal R=\partial_tU-\nu\Delta_hU+U\cdot\nabla U,
 \qquad
 q_* =\int_0^T\|c_U(t)\|_{X_{2\rho_0,m}}\,dt,
 \qquad
 r_* =\int_0^T\|\mathcal R(t)\|_{X_{2\rho_0,m}}\,dt.
\]
Here $c_U$ is defined by \eqref{eq:analytic-stability-c}.
Assume that $U$, $\Lambda_hU$, $ZU$, $c_U$, and $\mathcal R$ are
continuous on $[0,T]$ with values in $X_{2\rho_0,m}$.
There are constants $a_*>0$, $r_0>0$, and
$C_{\mathrm{stab}}<\infty$, depending only on
$A,T,\nu,\rho_0,m$, with the following property:
if
\begin{equation*}
                       q_*\leq a_*\rho_0,
                 \qquad r_*\leq r_0,
\end{equation*}
then the zero-rotation horizontally viscous Navier--Stokes equation
with initial datum $U(0)$ has a smooth solution on $[0,T]$, and
\begin{equation}\label{eq:analytic-stability-conclusion}
 \sup_{0\leq t\leq T}\|u(t)-U(t)\|_{X_{\rho_0/2,m}}
                              \leq C_{\mathrm{stab}} r_*.
\end{equation}
In particular,
\[
 \sup_{0\leq t\leq T}\|u(t)-U(t)\|_\infty
                              \leq C C_{\mathrm{stab}}r_*.
\]
The constants depend only on $A,T,\nu,\rho_0,m$, and are independent
of the size of the finite norm in \eqref{eq:analytic-stability-v0}.
\end{proposition}

\begin{proof}
Constants denoted by $C$ depend on the fixed operator radius and $m$.

\paragraph{Regularization.}
For $0<\eta\leq1$, solve
\begin{equation}\label{eq:analytic-stability-reg}
 \partial_tu^\eta-\nu\Delta_hu^\eta
   +J_\eta\mathbb P\bigl(w^\eta\cdot
                                      \nabla w^\eta\bigr)=0,
 \qquad w^\eta=J_\eta u^\eta,
 \qquad u^\eta(0)=U(0).
\end{equation}
Lemma~\ref{lem:friedrichs-global-smooth} gives a global smooth
solution of \eqref{eq:analytic-stability-reg}, preserving divergence
and reflection parity.  It also belongs locally to $X_{\rho,m}$
whenever the initial datum does.

Let $e=u^\eta-U$, $\bar U=J_\eta U$, and
$\bar e=J_\eta e$.
Put $\mathcal N(f)=f\cdot\nabla f$ for divergence-free fields, and define
\[
 F_\eta(U)
   =\mathbb P\mathcal R+
     J_\eta\mathbb P \mathcal N(J_\eta U)-\mathbb P\mathcal N(U).
\]
The error equation has forcing $-F_\eta(U)$, because
$\partial_tU-\nu\Delta_hU+\mathbb P\mathcal N(U)=\mathbb P\mathcal R$.
For the fixed smooth profile $U$, strong convergence of the Gaussian regularization at a
slightly smaller radius gives
\begin{equation}\label{eq:analytic-stability-consistency}
 \int_0^T\|F_\eta(U)(t)\|_{X_{\rho_0,m}}\,dt
                                  \leq C r_*+o_\eta(1).
\end{equation}
This convergence follows from the algebra estimate, convergence of
$J_\eta U$ and its horizontal and conormal derivatives, and
boundedness of $\mathbb P$.
The larger analytic radius in \eqref{eq:analytic-stability-A} controls the
conormal derivative tails; $U$ is fixed in this limiting step.

\paragraph{The analytic estimate.}
Use a decreasing radius $\rho(t)$ with $\rho(0)=\rho_0$, to be chosen
below.  Write
\[
 E_j(t)=\|\Lambda_h^je(t)\|_{X_{\rho(t),m}},\quad j=0,1,2.
\]
Also set
\[
 q(t)=\|c_U(t)\|_{X_{2\rho_0,m}}.
\]
The nonlinear difference inside the outer operator
$J_\eta\mathbb P$ is
\begin{align*}
 \mathcal N(\bar U+\bar e)-\mathcal N(\bar U)
 ={}&\bar U_h\cdot\nabla_h\bar e
       +\bar e_h\cdot\nabla_h\bar U
       +\bar e_h\cdot\nabla_h\bar e\\
    &+c_{\bar U}Z\bar e+c_{\bar e}Z\bar U+c_{\bar e}Z\bar e.
\end{align*}
Estimates \eqref{eq:analytic-stability-c}--
\eqref{eq:analytic-stability-c-J} imply
\[
 \|c_{\bar U}\|_{X_\rho}\leq Cq(t),\qquad
 \|c_{\bar e}\|_{X_\rho}\leq CE_1.
\]
The operator, product, and $Y_{\rho(t),m}(e)$ estimates consequently give
\begin{equation}\label{eq:analytic-error-differential}
\begin{aligned}
 E_0'+\nu E_2
 \leq{}& C\bigl((A+1+E_0)E_1+(A+q(t))E_0
                         +\|F_\eta(U)\|_{X_\rho}\bigr)
 \\
       &+\bigl(\rho'+C(q(t)+E_1)\bigr)Y_{\rho(t),m}(e).
\end{aligned}
\end{equation}
The term $C(q+E_1)E_0$ produced by
$Y_{\rho(t),m}(J_\eta e)\leq C(Y_{\rho(t),m}(e)+E_0)$ is included in the first line.
For completeness, the positive damping $\nu E_2$ follows by taking
the vertical supremum at each horizontal frequency: horizontal heat
then acts as the scalar factor $e^{-\nu t|\xi_h|^2}$.
After the factorial sum, the derivative of the moving weight supplies
exactly $\rho'Y_{\rho(t),m}(e)$.  The inequality can be read as an upper Dini
derivative inequality and integrated; it requires no differentiability
of an $L^\infty$ norm.

Choose $C_0$ larger than the constants multiplying the conormal
derivative norm in \eqref{eq:analytic-error-differential} and in
the vertical-derivative estimate leading to
\eqref{eq:vertical-derivative-differential}, and put
\begin{equation}\label{eq:analytic-stability-radius}
 \rho'=-\frac{\rho_0}{16T}-C_0(q(t)+E_1(t)).
\end{equation}
The constant term gives a strict decrease of the radius and
contributes $\rho_0/16$ to its loss on $[0,T]$.
On the bootstrap interval $E_0\leq1$, use
$E_1\leq\sqrt{E_0E_2}$ and Young's inequality to obtain
\begin{equation}\label{eq:analytic-stability-closed}
 E_0'+\frac\nu2E_2+\frac{\rho_0}{16T}Y_{\rho(t),m}(e)
 \leq C\left(1+A+\frac{(1+A)^2}{\nu}\right)E_0
             +C\|F_\eta(U)\|_{X_{\rho_0,m}}.
\end{equation}
Here $q(t)\leq CA$ follows already from
\eqref{eq:analytic-stability-c} and
\eqref{eq:analytic-stability-A}.

Since $e(0)=0$, Gronwall's inequality and integration of
\eqref{eq:analytic-stability-closed} give
\begin{equation}\label{eq:analytic-stability-quantitative}
 \sup_{t\le T}E_0(t)+\int_0^TE_2(t)\,dt+\int_0^TE_1(t)\,dt
 \le C_{A,T,\nu,\rho_0,m}
            \int_0^T\|F_\eta(U)(t)\|_{X_{\rho_0,m}}\,dt.
\end{equation}
For the last term, use
\[
 \int_0^TE_1(t)\,dt
 \le\sqrt{T\left(\sup_{t\le T}E_0(t)\right)\int_0^TE_2(t)\,dt}.
\]
Choose $a_*$ and $r_0$ sufficiently small.  By
\eqref{eq:analytic-stability-consistency}, for all sufficiently
small $\eta$, estimate \eqref{eq:analytic-stability-quantitative}
then implies
\[
 \sup_{t\le T}E_0(t)<\frac12,\qquad
 C_0q_*+C_0\int_0^TE_1(t)\,dt<\frac{\rho_0}{4}.
\]
Equation \eqref{eq:analytic-stability-radius} consequently gives
$\rho(t)>\rho_0/2$ on $[0,T]$.  Both bootstrap bounds improve,
with constants depending only on $A,T,\nu,\rho_0,m$.

\paragraph{Continuation of the analytic norm.}
For fixed $\eta$, the regularized nonlinearity is locally
Lipschitz in each fixed $X_\rho$.
Suppose the analytic interval ends at $0<T_*<T$.  Since $\rho$
is decreasing, extend it to $T_*$ by
$\rho(T_*)=\lim_{t\uparrow T_*}\rho(t)$.
Estimates \eqref{eq:analytic-stability-radius}--\eqref{eq:analytic-stability-quantitative}
give $\rho(T_*)>\rho_0/2$ and bound $u^\eta(t)$ in
$X_{\rho(T_*),m}$ for $t<T_*$, because $\rho(t)\ge\rho(T_*)$.
The regularized forcing is therefore bounded in that fixed Banach
space.  Duhamel's formula for the horizontal heat semigroup, which is
strongly continuous on this Wiener space, gives a strong norm limit
at $T_*$.  Its integrable $\Lambda_h$ smoothing bound
$C(t-s)^{-1/2}$ also gives a finite horizontal first-derivative norm
at this positive endpoint, so the radius equation can be restarted.
Restart the fixed-radius regularized equation at that limit.
Since $\rho'\leq-\rho_0/(16T)<0$, every subsequent time $t>T_*$
has a strictly smaller radius than $\rho(T_*)$, so the conormal derivative sums are
finite and the same integrated estimate applies again.  This rules
out termination of the analytic interval before $T$.
One can perform this argument first with a strictly smaller initial
radius and then pass upward to $\rho_0$ by monotone convergence.

\paragraph{Ordinary vertical derivatives.}
Put $v=\partial_z u^\eta$, so
$\partial_z w^\eta=J_\eta v$.
Differentiating \eqref{eq:analytic-stability-reg} and using
$\partial_z w^\eta_3=-\operatorname{div}_h w^\eta_h$
gives the exact equation
\begin{align*}
 \partial_tv-\nu\Delta_hv+J_\eta\mathbb P\bigl(
   &w^\eta_h\cdot\nabla_hJ_\eta v
    +c_{w^\eta}ZJ_\eta v
 \\
   &+(J_\eta v)_h\cdot\nabla_hw^\eta
     -(\operatorname{div}_hw^\eta_h)J_\eta v
                                      \bigr)=0.
\end{align*}
In particular, the equation is linear in $v$.
At the same radius as in \eqref{eq:analytic-stability-radius}, put
$V_j=\|\Lambda_h^jv\|_{X_{\rho(t),m}}$ for $j=0,1,2$.
On the bootstrap interval $E_0\le1$, the operator and product
estimates give, uniformly in $\eta$,
\begin{align*}
 \|w_h^\eta\cdot\nabla_hJ_\eta v\|_{X_\rho}
   &\le C(A+1)V_1,\\
 \|c_{w^\eta}ZJ_\eta v\|_{X_\rho}
   &\le C(q(t)+E_1)\bigl(Y_\rho(v)+V_0\bigr),\\
 \|(J_\eta v)_h\cdot\nabla_hw^\eta
       -(\operatorname{div}_hw_h^\eta)J_\eta v\|_{X_\rho}
   &\le C(A+E_1)V_0.
\end{align*}
Indeed, $\|w^\eta\|_{X_\rho}\le C(A+E_0)$,
$\|\Lambda_hw^\eta\|_{X_\rho}\le C(A+E_1)$, and
\eqref{eq:analytic-stability-c-J} gives
$\|c_{w^\eta}\|_{X_\rho}\le C(q(t)+E_1)$.
Thus the coefficient of $Y_\rho(v)$ contains only $q(t)+E_1$;
the remaining terms are linear in $V_0$ or $V_1$.
Boundedness of the outer operator $J_\eta\mathbb P$ and the moving
radius identity yield
\[
 V_0'+\nu V_2
 \le C(A+1)V_1+C(A+q(t)+E_1)V_0
      +\bigl(\rho'+C(q(t)+E_1)\bigr)Y_\rho(v).
\]
The choice \eqref{eq:analytic-stability-radius} absorbs the last
term at the same radius.  Applying
$V_1\le\sqrt{V_0V_2}$ and Young's inequality gives
\begin{equation}\label{eq:vertical-derivative-differential}
 V_0'+\frac\nu2V_2
 \leq C\left(\frac{(A+1)^2}{\nu}+A+1+q(t)+E_1(t)\right)V_0.
\end{equation}
Consequently
\begin{equation}\label{eq:analytic-stability-V-bound}
 \sup_{0\leq t\leq T}V_0(t)
    \leq C_{A,T,\nu,\rho_0,m}
                              \|\partial_zU(0)\|_{X_{\rho_0,m}}.
\end{equation}
The constant is independent of $\eta$ and of the size of the
initial vertical derivative.  The estimate can first be proved
at radius $\rho(t)-\delta$, then passed to the displayed radius
by letting $\delta\downarrow0$.

It follows from \eqref{eq:analytic-stability-quantitative} and
\eqref{eq:analytic-stability-V-bound} that
\begin{equation*}
 \int_0^T\|\nabla w^\eta(t)\|_\infty\,dt<\infty
\end{equation*}
with a bound uniform in $\eta$ for this fixed initial datum.
Self-adjointness of $J_\eta$ and its commutation with ordinary
derivatives give, for every fixed integer $k\geq3$,
\[
 \frac{d}{dt}\|u^\eta\|_{H^k}^2
       +2\nu\|\nabla_hu^\eta\|_{H^k}^2
 \leq C_k\|\nabla w^\eta\|_\infty
                                      \|u^\eta\|_{H^k}^2.
\]
Indeed, the nonlinear pairing becomes the usual transport pairing
of $w^\eta$ with itself, and
$\|w^\eta\|_{H^k}\leq\|u^\eta\|_{H^k}$.
Thus the ordinary Sobolev norms needed to pass to the limit
are bounded uniformly in $\eta$ on $[0,T]$.
Their dependence on the fixed ordinary Sobolev norms is confined
to this continuation step.

\paragraph{Passage to the limit.}
The regularized equation may be viewed as the exact equation with
the added forcing
\[
 \mathcal S_\eta
   =J_\eta\mathbb P\mathcal N(J_\eta u^\eta)
                                  -\mathbb P\mathcal N(u^\eta).
\]
For every fixed integer $k\geq5$, the Gaussian multiplier inequality
$\|(J_\eta-I)g\|_{H^{s-1}}
\leq C\sqrt\eta\|g\|_{H^s}$ and Sobolev products give
\[
 \|\mathcal S_\eta\|_{H^{k-2}}
                     \leq C_k\sqrt\eta
                                      \|u^\eta\|_{H^k}^2.
\]
For example, split this difference into
$(J_\eta-I)\mathbb P\mathcal N(u^\eta)$ and
$J_\eta\mathbb P
 [\mathcal N(J_\eta u^\eta)-\mathcal N(u^\eta)]$.
The first uses $\mathcal N(u^\eta)\in H^{k-1}$, and the second uses
$J_\eta u^\eta-u^\eta=O(\sqrt\eta)$
in $H^{k-1}$.
The usual smooth difference energy estimate in $H^{k-2}$, applied
to $u^\eta-u^{\eta'}$, therefore gives
\[
 \sup_{0\leq t\leq T}
 \|u^\eta(t)-u^{\eta'}(t)\|_{H^{k-2}}
                \leq C_{U,k,T}(\sqrt\eta+\sqrt{\eta'}).
\]
The constant is finite for the fixed datum by the uniform $H^k$ bound.
Thus the solutions converge strongly in $C([0,T];H^{k-2})$ for
every fixed $k\geq5$ to a smooth solution $u$ of the original
equation with initial datum $U(0)$.  In particular the convergence
is uniform in space and time, so the pointwise error estimate
passes to the limit.

To pass the analytic norm to the limit, fix a finite
conormal sum $n\leq J$, a horizontal frequency ball $|\xi_h|\leq B$,
and a vertical interval $|z|\leq R$.
Strong convergence in $C_tH^{J+2}$ implies convergence of each
$Z^n\mathcal F_h(u^\eta-u)$ in
$C_tL^1_{|\xi_h|\leq B}L^\infty_{|z|\leq R}$: use Cauchy--Schwarz
on the frequency ball, the one-dimensional vertical Sobolev
inequality, and boundedness of the polynomial coefficients of $Z^n$
on the fixed vertical interval.
The truncated analytic norm therefore satisfies the same
bound in the limit.  Let $J$, $B$, and $R$ tend to infinity and use monotone
convergence of the positive norm expressions.
Passing to the limit in \eqref{eq:analytic-stability-quantitative}
and using \eqref{eq:analytic-stability-consistency} proves
\eqref{eq:analytic-stability-conclusion}.
\end{proof}

\subsection{Construction of the approximate solution}

Set
\[
 G(z)=(1-2z^2)e^{-z^2},\qquad
 h(z)=\pi^{-1/2}e^{-z^2}.
\]
The profile $G$ has zero integral and $G(0)=1$, while
$\int h=1$.  Both functions are bounded under complex
dilations with $|\arg z|<\pi/4$.  Since
$\widehat G(\zeta)=c\zeta^2e^{-\zeta^2/4}$, a direct calculation
gives, for $j\le k$,
$|\langle G(2^j\cdot),G(2^k\cdot)\rangle_{\dot H^{1/2}}|
\le C2^{3(j-k)}$.  Consequently
\[
 F_N(z)=N^{-1/2}\sum_{j=1}^N G(2^jz)
\]
satisfies
\begin{equation*}
 \|F_N\|_{\dot H^{1/2}}\le C,
 \quad \|F_N\|_2\le C N^{-1/2},
 \quad F_N(0)=\sqrt N.
\end{equation*}
For each fixed integer $k\ge1$, one also has
\begin{equation*}
 \|Z^kF_N\|_\infty\le C_kN^{-1/2},
 \qquad \|Z^kG(M\cdot)\|_\infty\le C_k.
\end{equation*}
Indeed, $Z^kG$ vanishes at zero and decays at infinity, so its dyadic
dilation sum is uniformly bounded. For the Gaussian difference the
constants have factorial bounds on any strictly smaller complex
dilation sector. Evenness and zero integral imply that the primitives
from zero of both profiles are odd and decay at both infinities.

Let $K>0$, $L\ge R\ge1$, and $N,M\ge1$.  For $\delta>0$, put $b=\delta/\sqrt{LR}$.
Then $(K/\sqrt N)F_N(0)=K$.

\paragraph{The scalar equation and the divergence correction.}

Choose real cutoffs $\chi,\psi\in C_c^\infty(\mathbb R)$ nonnegative and equal to one
on $[-4,4]$. Let
\[
 \Psi_0(x_h)=\chi(x_1/L)\psi(x_2/R),\qquad
 A_0(x_h)=\chi(x_1/L)\psi(x_2/R)\sin x_2,
 \qquad \Psi(t)=e^{\nu t\Delta_h}\Psi_0.
\]
For each real parameter $\vartheta$ in a fixed bounded interval containing
the range of $G$, solve the scalar linear parabolic equation
\begin{equation*}
 \partial_t V+b\vartheta \Psi(t,x_h)\partial_2V-\nu\Delta_hV=0,
 \qquad V(0,x_h,\vartheta)=A_0(x_h).
\end{equation*}
Define the preliminary horizontal components
\[
 v_1(t,x_h,z)=\frac{K}{\sqrt N}F_N(z)V(t,x_h,G(Mz)),
 \qquad
 v_2(t,x_h,z)=b\Psi(t,x_h)G(Mz).
\]
They satisfy the exact identities
\begin{equation*}
 \partial_tv_1+v_2\partial_2v_1-\nu\Delta_hv_1=0,
 \qquad \partial_tv_2-\nu\Delta_hv_2=0.
\end{equation*}
The vertical mean of $v_1$ is generally nonzero for $t>0$, despite
$\int F_N=\int G=0$. Set
\[
 I(t,x_h)=\int_{\mathbb R}v_1(t,x_h,z)\,dz.
\]
Using the Gaussian $h$ fixed above, set
\begin{equation*}
 U_1=v_1-Ih,
 \qquad U_2=v_2,
 \qquad
 U_3(t,x_h,z)=-\int_0^z
               (\partial_1U_1+\partial_2U_2)(t,x_h,y)\,dy.
\end{equation*}
Then $U$ is divergence-free, $U_h$ is even in $z$, and $U_3$ is odd.
Both horizontal components have zero vertical mean. Hence $U_3$
decays at both vertical infinities and $U$ has finite energy.
At $t=0$, $I=0$.  Thus $U(0)$ has the prescribed horizontal
components, and the vertical mean compensation preserves decay
of $U_3$ throughout the construction.

The approximate solution satisfies the following bounds, uniformly in
the vertical concentration parameters.

\begin{lemma}[Analytic bounds for the approximation]
\label{lem:analytic-localization-approximation}
Fix $T,\nu>0$, $K>0$, and an integer $m\ge0$.
After decreasing $\rho_*$ from Lemma~\ref{lem:endpoint-analytic-calculus}
if necessary,
for $L,R\ge1$, $0<b\le1$, and $bT\le1$,
\begin{align}
 \sup_{0\le t\le T}\left(
 \sum_{j=0}^2\|Z^jU\|_{X_{\rho_*,m}}
       +\|\nabla_hU\|_{X_{\rho_*,m}}\right)
       &\le C_{m,K,T,\nu},
 \notag\\
 \sup_{0\le t\le T}\|c_U\|_{X_{\rho_*,m}}
 &\le C_{m,K,T,\nu}\left(K/L+b/R+Kb/(ML)\right),
 \notag\\
 \sup_{0\le t\le T}\|\mathcal R\|_{X_{\rho_*,m}}
 &\le C_{m,K,T,\nu}\left(L^{-1}+bR^{-1}+M^{-1}\right),
 \label{eq:analytic-localization-residual}
 \end{align}
where $c_U=U_3/z$ and $\mathcal R=U_t+U\cdot\nabla U-\nu\Delta_hU$.
More precisely,
\[
 \|ZU_1\|_{X_{\rho_*,m}}
 \le C_{m,K,T,\nu}(K/N+KbT+M^{-1}),
 \qquad
 \|ZU_2\|_{X_{\rho_*,m}}\le C_{m,T,\nu}b.
\]
The transport bound also holds with $Kb/(ML)$ replaced by $M^{-1}$
under the stated assumptions.  The constants and radius do not depend
on the ordinary vertical frequencies $2^N$ and $M$.
\end{lemma}

\begin{proof}
Fix a sector angle $0<\sigma<\pi/4$.  For an analytic function on
the positive and negative dilation sectors, write
\[
 \|v\|_{\mathcal S_{\sigma,m}}
 =\int_{\mathbb R^2}(1+|\xi_h|)^m
       \sup_{|\varphi|\le\sigma,\ z\in\mathbb R}
       |\mathcal F_hv(\xi_h,e^{i\varphi}z)|\,d\xi_h.
\]
Analyticity in $w=\log|z|$ and the Cauchy estimate give
\begin{equation}\label{eq:sector-to-analytic-X}
 \|v\|_{X_{\rho,m}}
 \le \frac{1}{1-\rho/\sigma}\,
             \|v\|_{\mathcal S_{\sigma,m}},
 \qquad 0\le\rho<\sigma.
\end{equation}
Apply Cauchy's estimate on circles of radius strictly less than
$\sigma$, then let the radius tend to the sector angle.
The same assertion applied to $Zv$ controls the analytic norm of one
additional conormal derivative.  More generally, Cauchy's estimate on
a pair of strictly nested fixed sectors controls any fixed number of
additional conormal derivatives in a smaller analytic radius.

We first justify the uniform estimates for the scalar equation with a
complex parameter $\vartheta$.  On each fixed disk
$|\vartheta|\le H$, its solution $V(t,x_h,\vartheta)$ is holomorphic
in $\vartheta$ and satisfies
\begin{equation}\label{eq:scalar-Wiener-parameter}
 \sup_{t\le T,\,|\vartheta|\le H}
 \left(\|V\|_{\mathcal A_m}
       +L\|\partial_1V\|_{\mathcal A_m}
       +(bT)^{-1}\|\partial_\vartheta V\|_{\mathcal A_m}
 \right)\le C_{m,T,\nu,H}.
\end{equation}
The scalar Fourier equation gives
\begin{equation}\label{eq:scalar-Wiener-differential}
 \frac{d}{dt}\|V\|_{\mathcal A_m}+\nu\|\Lambda_h^2V\|_{\mathcal A_m}
 \le b|\vartheta|\,\|\Psi\|_{\mathcal A_m}\|\Lambda_hV\|_{\mathcal A_m}.
\end{equation}
Here the weight $(1+|\xi_h|)^m$ is submultiplicative.
Since $\|\Lambda_hV\|_{\mathcal A_m}\le \|V\|_{\mathcal A_m}^{1/2}\|\Lambda_h^2V\|_{\mathcal A_m}^{1/2}$ and
$\|\Psi(t)\|_{\mathcal A_m}\le \|\Psi_0\|_{\mathcal A_m}\le C_m$, Young's inequality yields
\[
 \frac{d}{dt}\|V\|_{\mathcal A_m}+\frac\nu2\|\Lambda_h^2V\|_{\mathcal A_m}
 \le C_{m,\nu}b^2|\vartheta|^2\|V\|_{\mathcal A_m}.
\]
The initial norms of $A_0$ and $\Psi_0$ are uniform in $L,R\ge1$.
Differentiating in $x_1$ produces the source
$-b\vartheta(\partial_1\Psi)\partial_2V$ and initial data of size
$O(L^{-1})$.  Differentiating in $\vartheta$ produces the source
$-b\Psi\partial_2V$ and zero initial data.  Estimate \eqref{eq:scalar-Wiener-differential}, at one additional
horizontal order where needed, applies to these forced equations and proves
\eqref{eq:scalar-Wiener-parameter}.  Further fixed horizontal and time
derivatives are bounded in the same way.  Holomorphy follows, for
example, by the convergent Volterra iteration of the heat formula on
short time intervals and continuation using
\eqref{eq:scalar-Wiener-differential}.

The mixed derivative estimate for the vertical primitive and the
mean correction retains the localization factor in the first horizontal direction.  Differentiating the scalar equation twice gives
\[
 (\partial_t-\nu\Delta_h+b\vartheta\Psi\partial_2)
                              \partial_1\partial_\vartheta V
 =-b\Psi\partial_2\partial_1V
   -b(\partial_1\Psi)\partial_2V
   -b\vartheta(\partial_1\Psi)
                              \partial_2\partial_\vartheta V,
 \qquad \partial_1\partial_\vartheta V(0)=0.
\]
Since $\|\partial_1\Psi\|_{\mathcal A_j}\le C_j/L$ at every fixed
horizontal order, \eqref{eq:scalar-Wiener-parameter} and the forced version of
\eqref{eq:scalar-Wiener-differential} imply
\begin{equation}\label{eq:scalar-mixed-parameter-localization}
 \sup_{|\vartheta|\le H}
      \|\partial_1\partial_\vartheta V(t)\|_{\mathcal A_m}
       \le C_{m,T,\nu,H}\frac{bt}{L},\qquad
 \sup_{t\le T,\,|\vartheta|\le H}
      \|\partial_t\partial_1\partial_\vartheta V(t)\|_{\mathcal A_m}
       \le C_{m,T,\nu,H}\frac bL.
\end{equation}
For the second bound, use the equation for
$\partial_1\partial_\vartheta V$ and two additional fixed
horizontal orders in \eqref{eq:scalar-Wiener-parameter}.  The same bounds with a changed constant
hold after any fixed horizontal derivatives.

For the composition with $G$, choose the parameter radius $H$ so that
\[
 2\sup_{|\varphi|\le\sigma,\,z\in\mathbb R}
       |G(e^{i\varphi}z)|\le H.
\]
The Banach-valued Cauchy estimate then gives
\[
 \sum_{n\ge0}\frac{(H/2)^n}{n!}
       \|\partial_\vartheta^nV(t,\cdot,0)\|_{\mathcal A_m}
 \le 2\sup_{|\vartheta|\le H}\|V(\vartheta)\|_{\mathcal A_m}.
\]
This estimate permits the supremum over the sector to be taken inside
the horizontal Fourier integral.  It also applies to the derivatives
of $V$ in \eqref{eq:scalar-Wiener-parameter} and
\eqref{eq:scalar-mixed-parameter-localization}, including any fixed
additional horizontal derivatives.

Gaussian decay in any closed smaller sector gives
\begin{align*}
 \sup_{\mathrm{sector}}|\frac{K}{\sqrt N}F_N|&\le CK,&
 \sup_{\mathrm{sector}}|\frac{K}{\sqrt N}ZF_N|&\le CK/N,\\
 \int_{\mathbb R}|\frac{K}{\sqrt N}F_N(e^{i\varphi}z)|\,dz&\le CK/N,&
 \sup_{\mathrm{sector}}(|G(Mz)|+|ZG(Mz)|)&\le C.
\end{align*}
The $ZF_N$ estimate follows because $ZG$ vanishes at zero and decays
at infinity, so its dyadic dilation sum is uniformly bounded.
The integral estimate follows by summing the factors $2^{-j}$.
The chain rule and \eqref{eq:scalar-Wiener-parameter} therefore yield
the sector bounds for $v_1,v_2$, their horizontal derivatives,
and $Zv_1,Zv_2$.

The vertical mean satisfies the following uniform bound.
Since $\int F_N=0$,
\[
 I=\frac{K}{\sqrt N}\int_{\mathbb R}F_N(z)
          \bigl[V(t,x_h,G(Mz))-V(t,x_h,0)\bigr]\,dz.
\]
The parameter derivative estimate, $|\frac{K}{\sqrt N}F_N|\le CK$, and
$\|G(M\cdot)\|_1=M^{-1}\|G\|_1$ imply
\begin{equation}\label{eq:analytic-compensator-simple}
 \|I\|_{\mathcal A_m}\le C_{m,T,\nu}KbT/M.
\end{equation}
The scalar equation for $\partial_\vartheta V$ and
\eqref{eq:scalar-Wiener-parameter}, at two additional horizontal
orders, give $\|\partial_t\partial_\vartheta V\|_{\mathcal A_m}
\le C_{m,T,\nu}b$ on the fixed parameter disk.
Combining this bound with
\eqref{eq:scalar-mixed-parameter-localization} in the integral for $I$
gives
\begin{equation}\label{eq:analytic-compensator-derivatives}
\begin{aligned}
 \sup_{t\le T}\bigl(\|I\|_{\mathcal A_{m+2}}
       +\|\partial_tI\|_{\mathcal A_m}\bigr)
       &\le C_{m,T,\nu}\frac{Kb}{M},\\
 \sup_{t\le T}\bigl(\|\partial_1I\|_{\mathcal A_{m+2}}
       +\|\partial_t\partial_1I\|_{\mathcal A_m}\bigr)
       &\le C_{m,T,\nu}\frac{Kb}{ML}.
\end{aligned}
\end{equation}
Here the factor $M^{-1}$ comes from the sector $L^1$ norm of
$G(M\cdot)$; the factor $L^{-1}$ is retained by the mixed parameter
estimate.  Multiplication by $h$ and $Zh$ is bounded in every fixed
smaller sector.

Contour rotation is justified by Gaussian decay and gives
\[
 \int_{\mathbb R}v_1(x_h,e^{i\varphi}s)\,ds
       =e^{-i\varphi}I(x_h),\qquad
 \int_{\mathbb R}h(e^{i\varphi}s)\,ds=e^{-i\varphi}.
\]
The rotated integral of $G(M\cdot)$ also vanishes.  Thus
$\int_{\mathbb R}U_a(x_h,e^{i\varphi}s)\,ds=0$ for $a=1,2$.
Evenness makes the two half-line integrals zero separately, so
\[
 U_3(x_h,e^{i\varphi}z)
   =e^{i\varphi}\int_z^\infty
               \operatorname{div}_hU_h(x_h,e^{i\varphi}s)\,ds,
 \qquad z>0.
\]
For the derivative bounds, the defining primitive can equivalently
be written as
\begin{equation}\label{eq:analytic-U3-primitive}
\begin{aligned}
 U_3(t,x_h,z)
 ={}&-\frac K{\sqrt N}\int_0^zF_N(y)
          \partial_1V(t,x_h,G(My))\,dy\\
    &+\partial_1I(t,x_h)\int_0^zh(y)\,dy
      -b\partial_2\Psi(t,x_h)\,ze^{-(Mz)^2}.
\end{aligned}
\end{equation}
The sector $L^1$ bound for $(K/\sqrt N)F_N$ is $CK/N$.
The scalar equation and \eqref{eq:scalar-Wiener-parameter} give
$\|\partial_1V\|_{\mathcal A_{m+2}}
 +\|\partial_t\partial_1V\|_{\mathcal A_m}\le C_{m,T,\nu}/L$;
the corresponding bounds for $\partial_2\Psi$ are $C_{m,T,\nu}/R$.
Consequently, \eqref{eq:analytic-U3-primitive} and
\eqref{eq:analytic-compensator-derivatives} imply
\begin{equation}\label{eq:analytic-U3-derivatives}
\begin{aligned}
 \sup_{t\le T}\bigl(&\|U_3\|_{\mathcal S_{\sigma,m}}
   +\|\partial_tU_3\|_{\mathcal S_{\sigma,m}}
   +\|\nabla_hU_3\|_{\mathcal S_{\sigma,m}}\\
   &+\|\Delta_hU_3\|_{\mathcal S_{\sigma,m}}
   +\|ZU_3\|_{\mathcal S_{\sigma,m}}\bigr)
 \le C_{m,T,\nu}
       \left(\frac K{NL}+\frac{Kb}{ML}+\frac b{RM}\right).
\end{aligned}
\end{equation}
For $ZU_3$, apply the preceding primitive bounds in a fixed larger
sector and use Cauchy's estimate in the smaller sector of angle
$\sigma$; no ordinary vertical derivative is taken.
Finally,
\[
 c_U(z)=-\int_0^1\operatorname{div}_hU_h(az)\,da
\]
and the bounds for the first horizontal derivative and the mean correction give
\begin{equation}\label{eq:analytic-normal-transport-bound}
 \sup_{t\le T}\|c_U\|_{\mathcal S_{\sigma,m}}
 \le C_{m,T,\nu}
       \left(\frac K L+\frac bR+\frac{Kb}{ML}\right).
\end{equation}
All sector angles are fixed independently of $L,R,N,M$.

The exact triangular equations for $v_1,v_2$ give the three
components of the unprojected residual:
\begin{equation}\label{eq:analytic-residual-components}
\begin{aligned}
 \mathcal R_1={}&v_1\partial_1v_1+c_UZv_1
   -h(\partial_t-\nu\Delta_h)I-Ih\partial_1v_1\\
   &-hU_h\cdot\nabla_hI-Ic_UZh,\\
 \mathcal R_2={}&v_1\partial_1v_2+v_2\partial_2v_2
                   +c_UZv_2-Ih\partial_1v_2,\\
 \mathcal R_3={}&\partial_tU_3-\nu\Delta_hU_3
       +U_h\cdot\nabla_hU_3+c_UZU_3.
\end{aligned}
\end{equation}
The scalar estimates also give
$\|\partial_1v_1\|_{\mathcal S_{\sigma,m}}\le C_{m,T,\nu}K/L$,
$\|\partial_1v_2\|_{\mathcal S_{\sigma,m}}\le C_{m,T,\nu}b/L$, and
$\|\partial_2v_2\|_{\mathcal S_{\sigma,m}}\le C_{m,T,\nu}b/R$.
The sector algebra estimate and
\eqref{eq:analytic-normal-transport-bound} therefore control the
principal terms of $\mathcal R_1,\mathcal R_2$ by
$C_{m,K,T,\nu}(L^{-1}+bR^{-1}+M^{-1})$.
Every term containing $I$ or a derivative of $I$ is bounded by
$C_{m,K,T,\nu}/M$ by
\eqref{eq:analytic-compensator-derivatives}.
For $\mathcal R_3$, use \eqref{eq:analytic-U3-derivatives} and the
uniform sector bounds for $U_h,c_U$.  Since $N,M,L,R\ge1$ and
$b\le1$, this yields
\[
 \sup_{t\le T}\sum_{a=1}^3
     \|\mathcal R_a\|_{\mathcal S_{\sigma,m}}
 \le C_{m,K,T,\nu}\left(L^{-1}+bR^{-1}+M^{-1}\right).
\]
Applying \eqref{eq:sector-to-analytic-X} with
$\rho_*\le\sigma/4$ proves
\eqref{eq:analytic-localization-residual}.

Finally, fix $L,R,N,M,b$.  Commuting the scalar equation with
horizontal polynomial weights and ordinary derivatives, and inducting
on their orders, gives continuity of $V$ and its fixed parameter and
time derivatives in the horizontal Schwartz topology, locally
uniformly on bounded parameter disks.  Composition with $G(Mz)$
and multiplication by $F_N(z)$ or $h(z)$ make $U_h$ Schwartz in all
variables.  The zero half-line means and the tail formula for $U_3$
give its weighted decay; ordinary derivatives follow from the
primitive formula \eqref{eq:analytic-U3-primitive}.  Hence $U$ and
its fixed time derivatives are continuous in the Schwartz topology.
These ordinary bounds may depend on the fixed concentration parameters.
\end{proof}

\subsection{The critical lower bound}

Let $W(t)=\partial_\vartheta V(t,\cdot,0)$.
Choose a fixed nonnegative, nonzero
$\beta\in C_c^\infty((-1,1))$ and define
\[
 \phi_{L,R}(x_h)
 =-\frac{\beta(x_1/L)\beta(x_2/R)\cos x_2}
          {\|\beta(x_1/L)\beta(x_2/R)\cos x_2\|_{L^2(\mathbb R^2)}}.
\]
For all sufficiently large $R$,
\begin{equation}\label{eq:endpoint-horizontal-test}
 \|\phi_{L,R}\|_2=1,
 \qquad \|\phi_{L,R}\|_1\le C\sqrt{LR},
 \qquad -\int\phi_{L,R}\cos x_2\,dx_h\ge c\sqrt{LR}.
\end{equation}
Indeed, the integrals of the weighted $\cos^2x_2$ factors are
comparable to $R$, by their mean $1/2$ and decay of the Fourier
transform of the fixed smooth weights.

\begin{lemma}[Lower bound for the first parameter derivative]
\label{lem:endpoint-scalar-core}
For fixed $T,\nu>0$, there are $R_0,c_T,C_T>0$, independent of
$L,R,b$, such that for $L\ge R\ge R_0$ and $0<b\le1$,
\begin{align}
 \int_{\mathbb R^2}\phi_{L,R}(x_h)W(T,x_h)\,dx_h
       &\ge c_TbT\sqrt{LR},
 \label{eq:endpoint-W-lower}\\
 \sup_{x_h}
 |V(T,x_h,\vartheta)-V(T,x_h,0)-\vartheta W(T,x_h)|
       &\le C_Tb^2T^2|\vartheta|^2.
 \label{eq:endpoint-V-remainder}
\end{align}
The second estimate holds for every
$|\vartheta|\le\|G\|_\infty$.
\end{lemma}

\begin{proof}
At parameter zero, $V(t,\cdot,0)=e^{\nu t\Delta_h}A_0$, and
\[
 W(T)=-b\int_0^T e^{\nu(T-t)\Delta_h}
                   \bigl(\Psi(t)\partial_2e^{\nu t\Delta_h}A_0\bigr)
                 \,dt.
\]
On the larger rectangle $|x_1|\le2L$, $|x_2|\le2R$, Gaussian
heat-kernel tails imply, uniformly for $0\le t\le T$,
\[
 \Psi(t)=1+O(e^{-cR^2/(\nu T)}),\qquad
 \partial_2e^{\nu t\Delta_h}A_0
      =e^{-\nu t}\cos x_2+O(e^{-cR^2/(\nu T)}).
\]
For the differentiated assertion one may place $\partial_2$ on
$A_0$ before applying heat: the difference
$\partial_2A_0-\cos x_2$ vanishes on the cutoff plateau and is
uniformly bounded.  All the displayed functions are uniformly
bounded on the entire plane.  Applying the remaining heat semigroup
to their difference and using the distance from the support of
$\phi_{L,R}$ to the complement of the larger rectangle gives
\[
 W(T,x_h)
 =-bT e^{-\nu T}\cos x_2
                  +O\bigl(bT e^{-cR^2/(\nu T)}\bigr)
 \quad\text{on }\operatorname{supp}\phi_{L,R}.
\]
Together with \eqref{eq:endpoint-horizontal-test}, this proves
\eqref{eq:endpoint-W-lower} after increasing $R_0$.

The parameter-differentiated scalar equations have sources
$-b\Psi\partial_2V$ and
$-2b\Psi\partial_2\partial_\vartheta V$, respectively.
The uniform horizontal Wiener estimates used in
Lemma~\ref{lem:analytic-localization-approximation}, at two additional
fixed horizontal orders, give
\[
 \|\partial_\vartheta V(t)\|_{W_h^{1,\infty}}\le C_Tbt,
 \qquad
 \|\partial_\vartheta^2V(T)\|_\infty\le C_Tb^2T^2
\]
on the fixed parameter interval.  Taylor's formula proves
\eqref{eq:endpoint-V-remainder}.
\end{proof}

Set $g=|D_z|G$, where $|D_z|$ is the one-dimensional Fourier
multiplier with symbol $|\zeta|$.  Since $G\ne0$, one has
$\|G\|_{\dot H^{1/2}(\mathbb R)}>0$.  Moreover, $g\in L^1$, $\int g=0$, and
$\int Gg=\|G\|_{\dot H^{1/2}(\mathbb R)}^2$.  For example,
$\widehat G(\zeta)=c\zeta^2e^{-\zeta^2/4}$ shows that
$g$ decreases at least as $O(|z|^{-4})$.
With $e_1=(1,0,0)$, define the vector-valued test function
\[
 \Phi_M(x_h,z)=e_1\phi_{L,R}(x_h)M g(Mz).
\]
Fourier change of variables gives
\begin{equation}\label{eq:endpoint-dual-uniform}
 \|\Phi_M\|_{\dot H^{-1/2}(\mathbb R^3)}\le \|G\|_{\dot H^{1/2}(\mathbb R)},
 \qquad
 |\langle e,\Phi_M\rangle|
             \le C\sqrt{LR}\|e\|_\infty.
\end{equation}
The first inequality in \eqref{eq:endpoint-dual-uniform} follows from
$((|\xi_h|/M)^2+|\zeta|^2)^{-1/2}\le|\zeta|^{-1}$ and
$\|\phi_{L,R}\|_2=1$; the second uses $g\in L^1$.

For fixed $K,T,L,R,N,b$, dominated convergence in $y=Mz$ yields
\begin{align*}
 \lim_{M\to\infty}\langle U(T),\Phi_M\rangle
 &=K\int\phi_{L,R}(x_h)
       \int g(y)\bigl[V(T,x_h,G(y))-V(T,x_h,0)\bigr]\,dy\,dx_h.
\end{align*}
Here $\frac{K}{\sqrt N}F_N(y/M)\to K$, $|\frac{K}{\sqrt N}F_N|\le CK$, and $g\in L^1$.
The subtraction of $V(T,x_h,0)$ is permitted by $\int g=0$.
The contribution of the vertical mean correction tends to zero by
\eqref{eq:analytic-compensator-simple}.
Estimates \eqref{eq:endpoint-W-lower}--\eqref{eq:endpoint-V-remainder}
therefore give, whenever $bT$ is sufficiently small depending on
$T,\nu$ and the fixed profiles,
\begin{equation}\label{eq:endpoint-pairing-signal}
 \liminf_{M\to\infty}\langle U(T),\Phi_M\rangle
        \ge c_1KbT\sqrt{LR}.
\end{equation}
Indeed the main term is $K\|G\|_{\dot H^{1/2}(\mathbb R)}^2\int\phi_{L,R}W(T)$.
The absolute value of the remainder is at most
\[
          C_TKb^2T^2\sqrt{LR}\int_{\mathbb R}|G(y)|^2|g(y)|\,dy.
\]

\begin{proof}[Proof of Theorem~\ref{thm:illposed-main} for $\alpha=1$, $s=1/2$, and $\Omega=0$]
Fix $\varepsilon,B,T_0>0$ and $0<T<\min(1,T_0)$.
We construct Schwartz data first, with strict inequalities for
the subsequent compact-support approximation.
Choose $0<\delta\le1$ sufficiently small that the contribution of the second horizontal component in
\eqref{eq:endpoint-input-small} is below $\varepsilon/8$.  Put
$b=\delta/\sqrt{LR}$ after $L,R$ are chosen.
Choose $K$ sufficiently large that
\begin{equation*}
                    c_1K\delta T>16\|G\|_{\dot H^{1/2}(\mathbb R)}B.
\end{equation*}
We apply the analytic estimates with $m=0$.
These choices fix the constant $A$ in
Lemma~\ref{lem:analytic-localization-approximation} and hence all
the constants in Proposition~\ref{prop:endpoint-analytic-stability}.
They are independent of the remaining parameters $L,R,N,M$.
Reduce the analytic radius in
Lemma~\ref{lem:analytic-localization-approximation}, if
necessary, so that $2\rho_0$ lies below the admissible radius in
Proposition~\ref{prop:endpoint-analytic-stability}.

Proposition~\ref{prop:endpoint-analytic-stability} and
\eqref{eq:analytic-localization-residual} imply, once the
proposition's smallness hypotheses are met,
\begin{equation}\label{eq:endpoint-exact-uniform-error}
 \sup_{0\le t\le T}\|u(t)-U(t)\|_\infty
       \le C_{K,T,\nu}\bigl(L^{-1}+bR^{-1}+M^{-1}\bigr).
\end{equation}
The initial vertical profiles are finite sums of Gaussian derivatives
and their decaying primitives.  Thus $\partial_zU(0)\in X_{\rho_0,0}$
for fixed $N,M$, after choosing $\rho_0$ below the profile radius.
The norm in \eqref{eq:analytic-stability-v0} is finite and does not
enter the constant in \eqref{eq:endpoint-exact-uniform-error}.
By \eqref{eq:endpoint-dual-uniform}, the error in the duality pairing is at
most
\begin{equation}\label{eq:endpoint-error-order}
 C_{K,T,\nu}
       \left(\sqrt{R/L}+\delta/R+\sqrt{LR}/M\right).
\end{equation}

First take $R$ large enough
for Lemma~\ref{lem:endpoint-scalar-core} and to make the middle term
of \eqref{eq:endpoint-error-order} small.  Next take $L\gg R$ to
make its first term small, to ensure the required smallness of $bT$,
and to satisfy the transport and residual smallness assumptions of
Proposition~\ref{prop:endpoint-analytic-stability}, except for the $M^{-1}$ terms.  This is possible
because
\[
 \int_0^T\|c_U\|_{X_{2\rho_0,0}}\,dt
       \le C_{K,T,\nu}(K/L+b/R+M^{-1}),
\]
and $\mathcal R$ satisfies
\eqref{eq:analytic-localization-residual}.

Choose $N$ only after $R,L$ have been fixed.  The initial field satisfies
\begin{equation}\label{eq:endpoint-input-small}
 \|U(0)\|_{H^{1/2}}
 \le C\left(
    \frac{K\sqrt{LR}}{\sqrt N}+\delta
       +\frac{K\sqrt{LR}}{LN}+\frac{\delta}{RM}\right).
\end{equation}
Since $\partial_z(ze^{-z^2})=G(z)$, the primitives from zero are
\[
 \int_0^zF_N(y)\,dy
 =N^{-1/2}\sum_{j=1}^N2^{-j}(2^jz)e^{-(2^jz)^2},
 \qquad
 \int_0^zG(My)\,dy=ze^{-(Mz)^2}.
\]
The one-dimensional $\dot H^{1/2}$ norm is invariant under
vertical dilation, and the $L^2$ norm gains the factor $2^{-j/2}$.
The geometric sums therefore give
\begin{align*}
 \left\|\int_0^zF_N(y)\,dy\right\|_2
 +\left\|\int_0^zF_N(y)\,dy\right\|_{\dot H^{1/2}}
 &\le C N^{-1/2},\\
 \left\|\int_0^zG(My)\,dy\right\|_2
 +\left\|\int_0^zG(My)\,dy\right\|_{\dot H^{1/2}}
 &\le C M^{-1}.
\end{align*}
Moreover
$\|\partial_1A_0\|_{H_h^{1/2}}+\|\partial_1A_0\|_2
 \le C\sqrt{LR}/L$ and
$\|\partial_2\Psi_0\|_{H_h^{1/2}}+\|\partial_2\Psi_0\|_2
 \le C\sqrt{LR}/R$.
The product bound
\[
 \|a(x_h)F(z)\|_{H^{1/2}(\mathbb R^3)}
 \le \|a\|_{H_h^{1/2}}\|F\|_2
             +\|a\|_2\|F\|_{\dot H_z^{1/2}}
\]
applied to \eqref{eq:analytic-U3-primitive} at $t=0$, where
$I(0)=0$ and $V(0)=A_0$, gives the last two terms in
\eqref{eq:endpoint-input-small}.  The horizontal components use the
same product bound, $\|F_N\|_{\dot H^{1/2}}\le C$,
$\|F_N\|_2\le C N^{-1/2}$, and scale invariance for $G(Mz)$.
Taking $N$ large makes the background terms less than
$\varepsilon/8$.  This choice does not change any of the preceding
analytic constants.  Finally, with this $N$ fixed, take finite $M$
large enough for the limit in \eqref{eq:endpoint-pairing-signal},
for the remaining smallness hypotheses, for the last term of
\eqref{eq:endpoint-error-order}, and for the last initial-norm term.
In particular $M$ may be arbitrarily large relative to $2^N$.

Make the three terms in \eqref{eq:endpoint-error-order} together
smaller than $\|G\|_{\dot H^{1/2}(\mathbb R)}B$ and make the finite-$M$ duality pairing
at least half its limiting lower bound.
Proposition~\ref{prop:endpoint-analytic-stability} supplies a smooth
exact solution on $[0,T]$
with initial datum $U(0)$, and
\[
 \|U(0)\|_{H^{1/2}}<\varepsilon/2,
 \qquad
 \langle u(T),\Phi_M\rangle> 4\|G\|_{\dot H^{1/2}(\mathbb R)}B.
\]
Duality and \eqref{eq:endpoint-dual-uniform} give
$\|u(T)\|_{\dot H^{1/2}}>4B$.

\paragraph{Compactly supported initial data.}
With all parameters fixed, choose even zero-mean functions
$F^{(j)},G^{(j)}\in C_c^\infty(\mathbb R)$ converging to $F_N$ and
$G(M\cdot)$, respectively, in the Schwartz topology.  Such an
approximation is obtained by an even cutoff followed by subtraction
of its small remaining integral times one fixed even compactly
supported function of integral one.

Define the approximating data by
\begin{align*}
 u_{0,1}^{(j)}&=\frac{K}{\sqrt N}A_0(x_h)F^{(j)}(z),\\
 u_{0,2}^{(j)}&=b\Psi_0(x_h)G^{(j)}(z),\\
 u_{0,3}^{(j)}&=-\frac{K}{\sqrt N}\partial_1A_0(x_h)
                         \int_0^zF^{(j)}(y)\,dy
               -b\partial_2\Psi_0(x_h)
                         \int_0^zG^{(j)}(y)\,dy.
\end{align*}
They are real, smooth, divergence-free, and compactly supported.
Evenness and zero integral make each displayed primitive zero
outside a compact interval.  On the even zero-mean Schwartz subspace,
the primitive from zero equals the decaying primitive and is
continuous in the Schwartz topology.  Indeed, on $z>0$ it also equals
$-\int_z^\infty f(y)\,dy$, with the analogous tail formula on the
negative half-axis; weighted decay seminorms follow from these tail
formulas, and derivatives follow directly from $f$.
Consequently
$u_0^{(j)}\to U(0)$ in every fixed $H^k$.

Let $u^{(j)}$ be the classical solution with initial datum
$u_0^{(j)}$, and put $w=u^{(j)}-u$.  At integer order $k\ge3$,
the classical Sobolev stability estimate gives, as long as
$\|w\|_{H^k}\le1$,
\[
 \frac{d}{dt}\|w\|_{H^k}
       \le C_k(1+\|u\|_{H^{k+1}})\|w\|_{H^k}.
\]
The highest transport term cancels by incompressibility and
horizontal dissipation is nonnegative.  Gronwall and the classical
continuation criterion extend the approximating solutions through
$T$ and show convergence in $C([0,T];H^k)$.
For all sufficiently large $j$,
$\|u_0^{(j)}\|_{H^{1/2}}<\varepsilon$ and
$\|u^{(j)}(T)\|_{\dot H^{1/2}}>3B$.

\end{proof}

\subsection{The prescribed Coriolis parameter}

\begin{proof}[Completion of Theorem~\ref{thm:illposed-main} at the endpoint]
Fix $\Omega\in\mathbb R$ and $\varepsilon,B,T_0>0$.
The zero-rotation construction in
Section~\ref{sec:endpoint-inflation} gives one compactly supported smooth
datum $f$ and a zero-rotation classical solution $v$ on $[0,T]$
such that
\[
              \|f\|_{H^{1/2}}<\varepsilon,\qquad
              \|v(T)\|_{\dot H^{1/2}}>B+2.
\]
For a small parameter $\omega$, let $v_\omega$ solve \eqref{eq:NSC}
with $\alpha=1$, rotation rate $\omega$, and the same datum.
The $H^3$ difference-energy and continuation argument preceding
\eqref{eq:negative-error}, with equal viscosities and parameter
forcing $-\omega\mathbb P(e_3\times v)$, yields
\[
 \sup_{t\le T}\|v_\omega(t)-v(t)\|_{H^3}
                       \le C_{v,T,\nu}|\omega|
\]
for all sufficiently small $|\omega|$.  In particular these
solutions exist through $T$.

Define
\[
 u_\lambda(t,x)
       =\lambda v_{\Omega/\lambda^2}(\lambda^2t,\lambda x).
\]
The time derivative, horizontal Laplacian, and convection term all
scale by $\lambda^3$, whereas the Coriolis multiplier has order
zero.  Thus $u_\lambda$ solves \eqref{eq:NSC} with the prescribed
$\Omega$ and the same $\nu$.  Its initial datum
$\lambda f(\lambda x)$ remains real, solenoidal,
and compactly supported.  For $\lambda\ge1$,
\begin{align*}
 \|\lambda f(\lambda\cdot)\|_{H^{1/2}}^2
   &=\int_{\mathbb R^3}(\lambda^{-2}+|\xi|^2)^{1/2}
                                  |\widehat f(\xi)|^2\,d\xi
     \le\|f\|_{H^{1/2}}^2,\\
 \|\lambda F(\lambda\cdot)\|_{\dot H^{1/2}}
   &=\|F\|_{\dot H^{1/2}},\qquad F\in\dot H^{1/2}.
\end{align*}
Consequently
\[
 \|u_\lambda(T/\lambda^2)\|_{\dot H^{1/2}}
 \ge\|v(T)\|_{\dot H^{1/2}}
       -C_{v,T,\nu}|\Omega|/\lambda^2>B
\]
for large $\lambda$.  Increasing $\lambda$ also ensures
$T/\lambda^2<T_0$.  Together with
Section~\ref{sec:subcritical-inflation}, this proves
Theorem~\ref{thm:illposed-main}.
\end{proof}

\section*{Acknowledgments}
The research of Zefu Feng was supported by the National Natural Science Foundation of China (Grant No. 12101095), the Science and Technology Research Program of Chongqing Municipal Education Commission (Grant No. KJQN202300542), and the Talent Program of Chongqing Normal University (Grant No. BWQB2023008).
Yinghui Wang was partially supported by the National Natural Science
Foundation of China (Grant No. 12401274) and the Natural Science
Foundation of Hunan Province (Grant No. 2024JJ6302).
Jiahong Wu was partially supported by the US National Science Foundation
(Grant Nos. DMS-2104682 and DMS-2309748).

\section*{Declarations}

\noindent\textbf{Data availability.} No datasets were generated or analyzed in
this mathematical study.

\medskip
\noindent\textbf{Conflict of interest.} The authors declare no conflict of
interest.

\end{document}